\documentclass[preprint,12pt]{elsarticle}

 \usepackage{epsfig}

\usepackage{color}

\usepackage{amssymb}

\DeclareSymbolFont{msbm}{U}{msb}{m}{n}
\DeclareMathSymbol{\R}{\mathalpha}{msbm}{'122}

\begin{document}
\begin{frontmatter}
 


\title{ Simulation of a moving liquid droplet inside a rarefied gas region  \tnoteref{thanks}}
\tnotetext[thanks]{This work was partially supported by the German Research 
Foundation (DFG), KL 1105/17-1, HA 2696/16-1}


\author[TUKL]{Sudarshan Tiwari\corref{cor1}}
\ead{tiwari@mathematik.uni-kl.de} 
\cortext[cor1]{Corresponding Author}
\author[TUKL]{Axel Klar}
\ead{klar@mathematik.uni-kl.de}
\author[TUD]{Steffen Hardt}
\ead{hardt@csi.tu-darmstadt.de}
\author[TUD]{Alexander Donkov}
\ead{aadonkov@csi.tu-darmstadt.de}

\address[TUKL]{Fachbereich Mathematik, TU Kaiserslautern,
         Gottlieb-Daimler-Strasse, 
         67663 Kaiserslautern, Germany}
\address[TUD]{Center of Smart Interfaces, TU Darmstadt, 
Petersenstr. 32, 64287, TU Darmstadt Germany}

\begin{abstract}
We study the dynamics of a liquid droplet inside a gas over a large range of the 
Knudsen numbers.  The moving liquid droplet is modeled by  
the incompressible Navier-Stokes equations, the surrounding rarefied gas by 
the Boltzmann equation. The interface boundary conditions between the 
gas and liquid phases are derived. The incompressible Navier-Stokes equations 
are solved by a meshfree 
Lagrangian particle method called Finite Pointset Method (FPM), 
and the Boltzmann equation by a DSMC 
type of particle method. To validiate the coupled solutions of  
the Boltzmann and the incompressible Navier-Stokes equations 
we have further solved the compressible and the incompressible Navier-Stokes 
equations in the gas and liquid phases, respectively. 
In the latter case both the compressible and the incompressible 
Navier-Stokes equations are also solved by the FPM. 
In the continuum regime the coupled solutions obtained from the 
Boltzmann and the incompressible Navier-Stokes equations match with the 
solutions obtained from the compressible and the incompressible 
Navier-Stokes equations. In this paper, we 
presented solutions in one-dimensional physical space. 
\end{abstract}


\end{frontmatter}


\section{Introduction}
\label{secintro}
 
Gas-liquid flows in small scale geometries have received considerable 
attention in the past few years, especially due to the rapid developments 
in micro- and nanofluidics. Recently, experiments have been performed in 
which a segmented flow of gas 
and liquid is studied in nanochannels \cite{KKH, OFdeBM,  PNYJetal}. 
The Knudsen number of these flows, i.e. the ratio of the mean free path of the particles and the 
characteristic length, is such that the Boltzmann equation 
needs to be solved to describe the transport processes in the gas phase, 
while usually the incompressible Navier-Stokes equations are sufficient 
to model the liquid phase. While numerical methods for solving the 
Boltzmann equation are well established \cite{Bird, NS95}, no efficient 
methods seem to exist to study the flow of a gas at rarefied conditions 
coupled to the liquid flow described by the incompressible Navier-Stokes 
equations. The aim of this 
article is to take first steps into that direction. 

We consider a liquid droplet inside a gas. A liquid droplet may move, 
for example, due to a pressure difference, so that we need to solve 
a moving interface problem. One may choose different methods to compute 
these types of two-phase flows with moving interface, however, 
meshfree Lagrangian particle methods 
seem to be one of the preferred choices for such problems. 
For the rarefied gas phase we solve the Boltzmann equation by a 
DSMC type of particle method \cite{NS95, BI89}. In DSMC type of 
methods the gas particles are Lagrangian particles and move with their 
molecular velocities. However, these methods are mesh-based since  
one divides a computational domain into cells and 
performs intermolecular collisions of particles inside cells. Moreover, 
macroscopic quantities like fluid density, velocity, temperature, etc. are stored in the cell centers.  
On the other hand, we solve the incompressible 
Navier-Stokes equations by the 
Finite Pointset Method(FPM) \cite{TK07, TKH09}, which is a 
meshfree Lagrangian particle method and is similar in character as Smoothed 
Particle Hydrodynamics (SPH) \cite{GM77}. Within the FPM we  
approximate a liquid domain by Lagrangian particles which move with the local 
fluid velocity. We note that in this article we utilize two types of particles 
which may confuse the readers. Therefore, we call the DSMC particles "gas 
molecules" and the FPM particles "gas- and liquid particles" for 
the gas- and liquid phases, respectively. 
FPM particles are numerical grid points that move with the 
fluid velocity and carry all fluid dynamic quantities like densities, 
velocities, pressures, etc. along with them. To solve this type of two-phase flow 
problem one has to decompose the computational domain into a liquid and a 
gas domain. The domain decomposition is performed by first 
generating the entire domain by regular cells which are DSMC cells. Then we 
generate liquid particles overlapping the DSMC cells. The  
interface between the two domains is determined from the liquid particles. 
In the one-dimensional case, it is easily determined by identifying the extreme ends of 
the liquid particle positions. For higher-dimensional problems the interface 
is determined by identifying the free surface particles among the 
liquid particles, see \cite{TK02, TK07} for details. 

To the best of the author's knowledge, this is the first attempt to couple numerically 
the incompressible Navier-Stokes equations and the Boltzmann equation 
in the rarefied regime. 
To validate the coupled solution of the equations we have 
further simulated the gas and liquid phases by solving the compressible and 
the incompressible Navier-Stokes equations, respectively. 
In the latter case the interface 
conditions are quite standard which have been reported in a number of articles 
in the past few years \cite{BKZ, HW, HN, KP, TK07}, see also the references therein. 
We have considered the same examples as 
presented in \cite{CFA00, FMO98}, where the flow has been 
modeled by the Euler or the Navier-Stokes equations.  
Since we solve both the compressible and the incompressible equations by a particle method, 
we use different flags to distinguish the gas and liquid particles. 
Interfaces between the fluids can be determined based on their 
flags or colors \cite{TK02, MORRIS}. 
The flags are defined initially, and particles carry them along and leave them unchanged during the simulations. 

As will be shown in the article, for small Knudsen numbers 
the solutions obtained from the coupling of the Boltzmann and the 
incompressible Navier-Stokes equations are very close to the 
solutions obtained from the coupling of  
the compressible and the incompressible Navier-Stokes equations. 
However, the same is no longer true for the larger Knudsen numbers. We 
present test cases with smaller as well as larger Knudsen numbers. 

The paper is organized as follows. In section 2, we introduce the 
mathematical model for the gas and the liquid phases. In addtion to that  
we derive the interface boundary conditions. 
In section 3, we describe the particle methods for solving the Boltzmann and 
the Navier-Stokes equations. Section 4 is devoted to the coupling algorithms for 
the gas and the liquid phases. The numerical tests are presented in 
section 5, and some concluding remarks are given in section 6.

\section{Mathematical Model}
\label{mathmodel}
\subsection{Gas phase: the Boltzmann and the compressible Navier-Stokes equations}

The Boltzmann equation describes  the time evolution of a distribution 
function $f(t,{\bf x},{\bf v})$ for particles of velocity ${\bf v}\in \R^3$ 
at point ${\bf x}\in \Omega  \subset \R^s, ( s = 1,2,3)$ and time $t\in \R_+$. 
It is given in nondimensional form as 
\begin{eqnarray}
{\frac{\partial f}{\partial t}} + {\bf v} \cdot \nabla_x f = \frac{1}{\epsilon} J(f,f),
\label{BE}
\end{eqnarray}
where $\epsilon$ is the Knudsen number, $J(f,f)$ is 
the collision operator which is given for hard spheres by 
\begin{equation}
J(f, f) = \int_{\R^3} \int_{S^2} \beta(|{\bf v}-{\bf w}|, {\bf n}) [f({\bf v}^{'})f({\bf w}^{'}) - f({\bf v}) f({\bf w})] 
d{\bf n} d{\bf w}, 
\label{collisionoperator}
\end{equation}
where $S^2$ is the unit sphere 
in $\R^3$, ${\bf n} \in S^2$ is the unit vector in the impact direction, 
$\beta$ is the collision cross section, 
$ f({\bf v}^{'}) = f(t,{\bf x}, {\bf v}^{'})$ and analogously for 
$f({\bf v})$ etc. 
The pairs $({\bf v},{\bf w})$ and 
$({\bf v}^{'},{\bf w}^{'})$ are the pre- and post collisional velocities of two 
colliding particles, given by 
\begin{equation}
{\bf v}^{'}= {\bf v}-{\bf n} \left[{\bf n}\cdot ({\bf v}-{\bf w})\right], \quad
{\bf w}^{'} ={\bf w}+{\bf n} \left[{\bf n}\cdot ({\bf v}-{\bf w})\right] .
\end{equation}

The collision operator has five collisional invariants $\psi({\bf v}) = 1, {\bf v}, |{\bf v}|^2/2$ satisfying 
\begin{equation}
\int_{\R^3} \psi ({\bf v})J(f,f)d{\bf v} = 0. 
\label{conservationproperty}
\end{equation} 
In other words, $J(f,f)$ locally satisfies the conservation laws for mass, momentum and energy. 

The basic quantities of interest are the 
macroscopic ones, like the density $ \rho(t,{\bf x})$, mean velocity 
${\bf u} ={\bf u}(t,{\bf x})$ and the total energy $E=E(t,{\bf x})$, 
and are defined as 
\begin{eqnarray}
\rho =  \int_{\R^3} f(t, {\bf x}, {\bf v}) d{\bf v}, \quad \quad 
{\bf u} = \frac{1}{\rho}\int_{\R^3} {\bf v} f(t, {\bf x}, {\bf v}) d{\bf v} \\
E= \frac{1}{\rho} \int_{\R^3} \frac {|{\bf v}|^2}{2}  f(t, {\bf x}, {\bf v}) d{\bf v} = 
\frac{1}{2}|{\bf u}|^2 + e, 
\end{eqnarray}
where $e$ is the internal energy, defined by 
\begin{equation}
e = \frac{1}{\rho} \int_{\R^3} \frac {|{\bf v} - {\bf u}|^2 }{2}  f(t, {\bf x}, {\bf v}) d{\bf v}. 
\end{equation}
Moreover, the pressure tensor $\varphi$ and heat flux ${\bf q}$ are defined by 
\begin{eqnarray}
\label{stensor}
\varphi &=& \int_{\R^3} ({\bf v} - {\bf u})\otimes ({\bf v} - {\bf u})  f(t, {\bf x}, {\bf v}) d{\bf v} \\
{\bf q} &=&  \int_{\R^3} \frac {|{\bf v} - {\bf u}|^2 }{2} ({\bf v} - {\bf u})  f(t, {\bf x}, {\bf v}) d{\bf v}. 
\label{hflux}
\end{eqnarray}

The gas pressure $p$ is defined as 
$ p = \frac{2}{3} \rho e$ for a monoatomic ideal gas. Furthermore, $ p = \rho R T$ 
holds, where $T$ is the temperature and $R$ is the gas constant. 
For more details  we refer
to \cite{CIP94, Sone07}. 

Multiplying the Boltzmann equation by its collisional invariants and then integrating 
with respect to ${\bf v}$ over $\R^3$ we obtain the following local conservation 
equations 
\begin{eqnarray}
 \label{momsys}
\frac{\partial \rho}{\partial  t} +\nabla \cdot (\rho {\bf u}) &=& 0\nonumber \\
\frac{ \partial (\rho  {\bf u})}{ \partial t} + \nabla \cdot (\rho {\bf u}\otimes {\bf u} +\varphi ) &=& 0 \\
\frac{\partial (\rho E) }{\partial t} + \nabla \cdot 
\left [ \rho E  {\bf u} 
+\varphi \cdot {\bf u} + {\bf q} \right ] &=& 0. \nonumber
\end{eqnarray}

For $\epsilon$ tending to zero, i. e. for small mean free paths, one 
can show that the phase space distribution function $f$ tends to 
the local Maxwellian \cite{Caf80}
\begin{eqnarray}
\label{maxwell}
f_{M}(t, {\bf x}, {\bf v})= \frac{\rho}
{(2\pi R T)^{3/2}} \;e^{- {\frac{|{\bf v} - {\bf u}|^{2}}{2 R  T}}}
\end{eqnarray}
and the system of moment equations (\ref{momsys}) tends to  the 
compressible Euler equations with the closure relations 
$\varphi = p I$ and ${\bf q} = {\bf 0}$. This can be verified from the 
asymptotic expansion of $f$ in $\epsilon$, where the zeroth order 
approximation gives the local Maxwellian distribution, 
and the first order approximation \cite{BGL} gives the 
Chapman-Enskog distribution
\begin{equation}
f_{CE}(t, {\bf x},{\bf v}) = f_M (t,{\bf x},{\bf v}) \left [ 1+\frac{2}{5} \frac{{\bf q}\cdot {\bf c}}{\rho(RT)^2} \left (
\frac {|{\bf c}|^2}{2RT}-\frac{5}{2} \right ) +  
\frac{1}{2}\frac{\tau:{\bf c}\otimes {\bf c}}{\rho (RT)^2} \right ], 
\end{equation}
where ${\bf c} = {\bf v}-{\bf u}$. At the same time, (\ref{momsys}) tends to the  
compressible Navier-Stokes equations with the closure relations 
\begin{equation}
\varphi = pI - \tau, \quad \quad  {\bf q} = - \kappa \nabla T, 
\label{momcloser}
\end{equation}
where 
\begin{equation}
\tau_{ij} = \mu \left (\frac{\partial u_i}{\partial x_j} + 
 \frac{\partial u_j}{\partial x_i} - \frac{2}{3} (\nabla\cdot {\bf u})~ \delta_{ij} 
\right ), \label{stresscloser}
\end{equation}
and $\mu = \mu(t,x)$ and $\kappa = \kappa(t,x)$ are the dynamic viscosity and thermal conductivity, respectively. They are of the order of $\epsilon$. 
For example, the first approximation for the viscosity and the thermal 
conductivity of a monatomic gas is given by \cite{VK65} 
\begin{equation} 
\mu = \frac{5}{16 d^2}\sqrt{\frac{mkT}{\pi}}, \quad \quad 
\kappa = \frac{15 k}{4m}\mu,
\label{mukappa}
\end{equation}
where $k$ is the Boltzmann constant, and $d$ and $m$ are the mass and the diameter of the molecules, respectively. 
In this paper we compute the parameters $\mu$ and $\kappa$ from the 
initial temperature, and use the corresponding values in the compressible Navier-Stokes equations. 

Since we solve the Navier-Stokes equations 
with a meshfree Lagrangian particle method, we re-express them in 
Lagrangian form with respect to the primitive variables as  
\begin{eqnarray}
\label{Lagrangian_NS}
\frac{d {\bf x}_g }{dt} &=& {\bf u}_g \nonumber \\
{\frac{d \rho_g}{d t}} &=& -\rho_g \nabla \cdot {\bf u}_g \nonumber \\
{\frac{ d {\bf u}_g}{d t}} &=&-\frac{1}{\rho_g} \nabla p_g + 
\nu_g \left [ \nabla^2{\bf u}_g + \frac{1}{3}\nabla(\nabla\cdot {\bf u}_g) \right ] 
\\
{\frac{ d T_g}{d t}} &=& \frac{1}{c_v\rho_g} \left[ -p_g ~ \nabla\cdot{\bf u}_g + 
(\tau_g\cdot\nabla)\cdot{\bf u}_g + {\kappa_g} \nabla^2 T_g \right ], \nonumber
\end{eqnarray}
where we have used the index $g$ for the gas quantities,  $\nu$ is the 
kinematic viscosity, $d/dt $ is the 
material derivative, $c_v$ is the specific heat 
at constant volume, given by $\frac{3}{2}R$ for a monoatomic gas. 
We note that we have expressed the internal energy of the gas as 
$e_g = c_v T_g$.

When introducing the specific heat into the energy equation an ideal gas was assumed. 
In addition to the system of equations 
(\ref{Lagrangian_NS}) we consider the equation of state 
\begin{equation}
p_g = \rho_g R T_g
\label{stateeqn}.
\end{equation}

\subsection{Liquid phase: Incompressible Navier-Stokes equations}

We consider incompressible flow inside the liquid phase. 
The governing equations can be obtained from the incompressible 
Navier-Stokes equations by 
assuming the liquid density $\rho_l$ to be constant. All liquid 
quantities are denoted with the index $l$. 
Moreover, one can express 
the internal energy $e_l$ 
for the liquid phase approximately by $e_l = c_p T_l$, where 
$c_p$ is the specific heat at constant pressure \cite{Lesieur}.  
We also solve the incompressible Navier-Stokes equations by a meshfree 
Lagrangian particle method. Therefore, we express these equations in 
Lagrangian form
\begin{eqnarray}
\label{incomp_NS}
\frac{d {\bf x}_l}{dt} &=& {\bf u}_l \nonumber \\
\nabla \cdot {\bf u}_l &=& 0 \nonumber \\
\frac{d {\bf u}_l}{d t} &=& -\frac{\nabla p_l}{\rho_l} + \nu_l \nabla^2 {\bf u_l} \\
{\frac{ d T_l}{d t}} &=& \frac{1}{c_p \rho_l}\left(\tau_l\cdot\nabla\right )\cdot{\bf u}_l + 
\frac{\kappa_l}{c_p \rho_l} \nabla^2 T_l,\nonumber
\end{eqnarray}
where $\tau_l$ is given by (\ref{stresscloser}) without the divergence of velocity term.
In many situations, the viscous dissipation term in the energy equation can be neglected. A detailed discussion about when it is justified to omit that term can be found in~\cite{MS07,Morini05}. In what follows, we limit ourselves to scenarios with negligible viscous heating, resulting in a simplified energy equation
\begin{equation}
\frac{d T_l}{d t} = \frac{\kappa_l}{c_p \rho_l} \nabla^2 T_l. 
\label{heateqn}
\end{equation}

\subsection{Initial and boundary conditions}
In this paper we consider a one-dimensional computational domain 
$\Omega = [a,b] \subset \R^1$, where 
$a$ is always zero and $b$ varies from $10^{-4}$ to $10^{-6}$.
The domain is initially decomposed into the gas domain $\Omega_g$ and 
liquid domain $\Omega_l = \Omega \setminus \Omega_g$. 
We consider cases where the liquid domain always remains inside of the full domain. 
Therefore, the boundaries $a$ and $b$ 
always belong to the gas domain. 
We prescribe boundary conditions for the gas at points $a$ and $b$. 
Moreover, there are interfaces between the liquid and the gas domains, 
and we have to further specify the interface boundary conditions, described in the next subsection. 

In the gas domain $\Omega_g$ we either solve the compressible Navier-Stokes 
equations or the Boltzmann equation. We assume that initially the gas is 
in thermal equilibrium with the values $\rho_g(0,x), u_g(0,x)$ 
and $T_g(0, x)$, which are the initial conditions 
for the compressible Navier-Stokes equations. 
If we solve the Boltzmann equation in 
$\Omega_g$, we prescribe the initial condition as a local Maxwellian with parameters 
$\rho_g(0,x), {\bf u}_g(0,x) = (u_g(0,x),0,0)$ and $T_g(0, x)$.
In $\Omega_l$ we solve the incompressible 
Navier-Stokes equations with initial conditions 
$p_l(0,x), u_l(0,x)$ and $T_l(0,x)$. 

\subsection{Interface boundary conditions}
To couple the liquid and gas phases one has to first determine the interface 
between two phases and then prescribe the interface boundary conditions. 
For solving the Boltzmann and the incompressible Navier-Stokes equations, 
we determine the interface as the free surface particles from the liquid domain.  
In the one dimensional case, we have two interfaces, which are 
given by the liquid particle position at the extreme left $(x_L)$ and the 
liquid particle position at the extreme right $(x_R)$. They can be tracked at 
every time step. 
When we solve the transport equations in both phases by the FPM, 
we assign different flags or colors for the particles in the compressible and 
incompressible phases. 
The particles of each phase 
carry the color function along with them, and the interface can be 
tracked with the help of the flags of the particles. 
We again determine the interface by identifying the liquid particles at the extreme left and right. 

Owing to the kinematic boundary condition at the interface, 
there is no penetration of particles from one phase to the other. This means that the convective terms for mass, 
momentum and energy transport are zero. Hence, all fluxes with the multiplicative 
factors ${\bf u}$ vanish.  Therefore, we have the following 
jump conditions for the momentum and energy fluxes 
the system (\ref{momsys}) 
\begin{eqnarray}
(\varphi_{11})_l &=& (\varphi_{11})_g \label{momjump} \\
(\varphi_{11} u + q)_l &=& (\varphi_{11}  u +  q)_g. \label{energyjump} 
\end{eqnarray}
Here we use a scalar quantity $q$ for the heat flux in the one-dimensional case. 
Moreover, we assume that the velocity and temperature at the interface are 
continuous, 
{\it i.e. }
\begin{eqnarray}
 u_l =  u_g. \label{velojump} 
\\
T_l = T_g. \label{tempjump}
\end{eqnarray}
Then from (\ref{momjump}), (\ref{energyjump}) and (\ref{velojump}) we have 
\begin{equation}
 q_l =  q_g \label{energyjump1}.
\end{equation}

\subsubsection{Interface boundary conditions for the compressible and the 
incompressible Navier-Stokes equations}
\label{interfaceNS_NS}

We use the closure relations (\ref{momcloser}) and (\ref{stresscloser}) in 
(\ref{momjump}) and (\ref{energyjump1}) and get the following continuity  
of the fluxes  
\begin{eqnarray}
\left ( p - \frac{4}{3}\mu \frac{\partial u }{\partial x} \right )_l &=& 
\left ( p - \frac{4}{3}\mu \frac{\partial u }{\partial x} \right )_g 
\label{compnsmomjump} \\ 
\left (
\kappa \frac{\partial T}{\partial x}\right )_l &=&
\left (
\kappa \frac{\partial T}{\partial x}\right )_g.
\label{compnshfluxjump}
\end{eqnarray}

Using the divergence-free condition for the liquid we obtain the 
interface boundary condition for the pressure as 
\begin{equation}
\left ( p  \right )_l =
\left ( p - \frac{4}{3}\mu \frac{\partial u }{\partial x} \right )_g.
\label{pressjump}
\end{equation}

We note that we use the condition (\ref{tempjump}) together with 
(\ref{compnshfluxjump}) for the thermodynamic equations. 

\subsubsection{Interface boundary conditions between the Boltzmann and the 
incompressible Navier-Stokes equations}

An important point to note in relation to the interface conditions between the Boltzmann and the Navier-Stokes domain is the fact that
while the quantities available on the gas side of the interface determine all of the quantities needed on the liquid side, the same is not true in the
opposite direction. Therefore, additional assumptions have to be made when computing the interfacial phase-space distribution in the gas from the density, velocity and temperature in the liquid. This requires a different treatment depending on whether information is passed from the gas to the liquid or vice versa.
As in the previous subsection \ref{interfaceNS_NS} we assume the continuity 
of the velocity and the temperature as given by (\ref{velojump}) and 
(\ref{tempjump}). 
In this case the closure relations (\ref{momcloser}) and 
(\ref{stresscloser}) are used only for the fluxes of the 
incompressible Navier-Stokes equations. 
Therefore, the continuity relations for the fluxes differ 
and are re-expressed in the form 
\begin{eqnarray}
(p)_l = \left(\varphi_{11}\right)_g \label{boltzpressjump} \\
-\left (  
\kappa \frac{\partial T}{\partial x}\right )_l = (q)_g, \label{boltzhfluxjump}
\end{eqnarray} 
where $(\varphi_{11})_g $ and $(q)_g$ are computed directly 
from the moments (\ref{stensor}) and (\ref{hflux}), respectively, 
of the solution of the Boltzmann equation. 
Hence, (\ref{tempjump}), (\ref{boltzpressjump}) and (\ref{boltzhfluxjump}) 
give the interface conditions from the gas into the liquid, where $T_g$ is  
computed from the solution of the Boltzmann equation. 

On the other hand, the interface boundary condition from the liquid into the 
gas is treated as follows. When gas molecules hit the interface, we apply 
the diffuse reflection condition with thermal accommodation, {\it i.e. } particles 
are reflected with the interface temperature and velocity into the 
gas domain, see the detailed descriptions in section \ref{BoltzNScouple}.

\section{Numerical methods}
\label{nummeth}

\subsection{Particle method for the Boltzmann equation}

For solving the Boltzmann equation we have used 
a variant of the DSMC method \cite{Bird}, developed in 
\cite{NS95, BI89}. The method is based on 
the time splitting of the Boltzmann equation. Introducing 
fractional steps one first solves the free transport equation (the 
collisionless Boltzmann equation) for one time step. During the 
free flow,  boundary and interface conditions are taken into account. 
In a second step (the collision step) the spatially 
homogenous Boltzmann equation 
without the transport term is solved. 
To solve the homogeneous Boltzmann equation the key point is to find 
an efficient particle approximation of the product distribution functions 
in the Boltzmann collision operator given only an approximation of the 
distribution function itself. To simulate this equation by a 
particle method an explicit Euler step is performed. 
To guarantee positivity of the distribution 
function during the collision step a restriction of the time step 
proportional to the Knudsen number is needed. That means that the 
method becomes exceedingly expensive for small Knudsen numbers.

\subsection{FPM for the compressible Navier-Stokes equations}
\label{fpm}
We solve the Navier-Stokes equations (\ref{Lagrangian_NS}) by the FPM. 
As already pointed out, the FPM is a 
meshfree Lagrangian particle method, 
where we approximate the spatial derivatives with the 
help of the weighted least squares method. In order to solve these 
equations by FPM, 
one first fills the computational 
domain by particles which can be considered as moving numerical grid points and 
then approximates the spatial derivatives occurring on  the right hand side  
of (\ref{Lagrangian_NS}) 
at each particle position from its neighboring  particles. 
This reduces the system of partial differential 
equations (\ref{Lagrangian_NS}) to a system of ordinary differential 
equations with respect to time. We solve the resulting ODE system 
with the help of a two-step Runge-Kutta method. The time steps for 
the compressible as well as the incompressible Navier-Stokes equations are 
restricted by the 
CFL condition and by the value of the transport coefficient  
$max\left [\mu_g/\rho_g, \kappa_g/(\rho_g c_v),  
\kappa_l/(\rho_l c_p) \right ]$. 
We refer to our earlier papers \cite{TKH09, Tiw01} for the details of 
the least squares approximation of the spatial derivatives. 
We note that we need to introduce a particle management scheme during the simulations. 
Because of the Lagrangian description of the method 
particles may accumulate or 
may thin out causing holes in the computational domain. This gives rise to   
some instabilities of the method. Therefore, we have to add or remove 
particles. In the one-dimensional case this task is quite simple. If the 
distance between a particle and its nearest neighbor is larger than 
$1.2$ times the 
initial spacing of particles, we add a new particle in the center. On the 
other hand, if two particles are closer than $0.2$ times the 
initial spacing we remove both of them and add a new particle at 
the mid point. However, these adding and remove factors may depend upon 
users. 
The fluid dynamic quantities 
of newly added particles are approximated from their neighbor particles with 
the help of the least squares method.

\subsection{FPM for the incompressible 
Navier-Stokes equations}

The incompressible Navier-Stokes equations (\ref{incomp_NS}) are 
also solved by the FPM. Since we consider one-dimensional physical space, 
the solution of the incompressible Navier-Stokes equations is simple. 
For higher-dimensional physical spaces we refer to \cite{TK02,TK07}. 

Taking the partial derivatives with respect to $x$ on both sides of 
the momentum equation (\ref{incomp_NS}) and using the divergence-free 
constraint, we obtain the one-dimensional Laplace equation for the pressure 
\begin{equation}
\frac{\partial ^2 p_l}{\partial x^2} = 0,
\end{equation} 
with the Dirichlet boundary conditions $p_L$ and $p_R$ on the left and right interfaces 
$x_L$ and $x_R$, respectively. 

Suppose we have $N$ liquid particles ${x_{l1},\cdots, x_{lN}}$ 
in $\Omega_l$. 
We have to solve the corresponding transport equations for every 
liquid particle. Let $x_{li}$ be the 
position of an arbitrary particle. 
The new pressure at time 
level $t^{n+1} = (n+1)\Delta t$, where $\Delta t$ is the time step, 
can be computed explicitely by 
\begin{equation}p^{n+1}_{li} = \frac{p_R^n - p_L^n}{x_R^n -x_L^n}x_{li}^{n} + \frac{p_R^n x_L^n - p_L^nx_R^n}{x_L^n - x_R^n}.
\label{incomppress}
\end{equation}

The new velocity is given by 
\begin{equation}
u^{n+1}_{li} = u^{n}_{li} - \frac {\Delta t}{\rho_l} (\frac{p_R^n-p_L^n}{x_R^n-x_L^n}).
\label{incompvelo} 
\end{equation}

The energy equation for the liquid is solved in the same way as in the 
case of the gas. 

We note that the equation for the velocity does not contain any contributions from the viscous term, which is due to the one-dimensional nature of the problem. The thermal 
diffusivity $\kappa_l/(\rho_l c_p)$ of the liquid is much smaller than  
that of the gas. 
Therefore, the time step for the gas phase 
is smaller than the one for the liquid. 
So, we use an explicit Euler step 
for the time integration to solve the energy equation of the liquid phase. 
However, the time steps are chosen as the smallest of the values obtained for the two phases. 

Finally, we compute the new 
liquid particle positions by
\begin{equation}
{x_{li}}^{n+1} = {x_{li}}^n + \Delta t {u}^n_{li}. \label{newposition}
\end{equation}

\section{Coupling of the gas and the liquid phase}

The model for the coupling of the two phases to be presented in the 
following contains some essential features of gas-surface interactions, but 
also relies on some simplifications. Essentially, we neglect evaporation and 
condensation phenomena, i.e. we assume mass conservation for both the gas and 
the liquid. While this could be an unsupported assumption for cases in which 
rarefaction is due to the reduced density of the gas phase, in the cases we 
study the reduction of the length scale down to micrometer dimensions causes 
rarefaction. Still, the results to be presentented in the following sections 
show that there exist local variations of the gas density and temperature. 
Therefore, strictly speaking, while the global (integrated) evaporation and 
condenstion mass fluxes could be very small, they may be of relevance 
locally. Nevertheless, two arguments show that even our comparatively 
simple model for gas-liquid interaction can make practically relevant 
predictions about gas-liquid flows. First of all, we study fast processes 
occuring on time scales below on microsecond where phase-change phenomena 
are usually negligible. Moreover, there exist non-volatile liquids for which 
even on longer time scales evaporation does not play a 
role \cite{Freemantle}. Therefore, with the necessary precations taken, the 
model presented in the following should provide a realistic picture of a 
specific class of gas-liquid flows.

\subsection{Coupling of the compressible and incompressible Navier-Stokes 
equations}
\label{NSNScouple}

We initially generate liquid and gas particles in the 
computational domain. We assign all initial values such as density, velocity, 
temperature, etc. to all particles. Moreover, we also assign different 
flags or colors to liquid and gas particles. For instance, we define a flag 
equal to 1 
for liquid particles and 2 for gas particles. These flags remain 
unchanged throughout the simulation. The gas equations are solved by  
gas particles and the liquid equations by liquid particles. 
We approximate the spatial derivatives at a gas particle from its 
neighbors and exclude all liquid particles from the neighbor list. This 
situation occurs near the interface. Similarly, we exclude the gas 
particles from the neighbor list of liquid particles. 
The interface particles at $x_L$ and $x_R$ are determined after 
every time steps via the leftmost and rightmost positions among the 
liquid particles.  
For higher dimensional cases, we need to spend 
more effort to determine the interface particles, but  
this is also quite straightforward, see \cite{TK02, TK07}. 

To couple both types of equations, we first solve the compressible 
Navier-Stokes equations for the gas particles. 
We consider the velocity $u_l$ and the temperature $T_l$ of the 
interface particles as Dirichlet boundary conditions.  
This gives 
the coupling conditions from the liquid phase into the gas phase. 

The coupling from the gas phase into the liquid phase is obtained 
as follows. 
After computing new quantities for the gas particles, we approximate the 
pressure and temperature at the interface positions $x_L$ and $x_R$. To  
approximate the pressure according to (\ref{pressjump}), we first extrapolate the 
pressure from the neighboring gas particles using the least squares method. 
Then we subtract the term containing the derivative of $u_g$ which is also 
computed from the neighboring gas 
particles from the least squares method. Subsequently we obtain the liquid pressure from (\ref{incomppress}) and the velocity from 
(\ref{incompvelo}) after determining the pressure on the interface. 

The approximation of the temperature at the interface is slightly 
different. At the interface the temperature and the heat flux are continuous  
according to (\ref{tempjump}) and (\ref{compnshfluxjump}), implying a jump condition for the temperature gradient.
We apply a two-sided interpolation satisfying these jump conditions as follows. 
Let $x_I$ be a interface particle 
with the temperature $T_I$. 
The particle $x_I$ has $m$ neighboring gas particles 
${x_{g1}, x_{g2}, \cdots, x_{gm}}$ and $n$ neighboring liquid particles 
${x_{l1}, x_{l2}, \cdots, x_{ln}}$. We write the Taylor expansions of 
$T_{li}$ and $T_{gi}$ around $T_I$ as 
\begin{eqnarray}
T_{gi} = T_I + dx_{gi} (\frac{\partial T_g}{\partial x_g}) + e_{gi}, \quad i = 1, \cdots, m \label{liqtaylor}  \\
T_{li} = T_I + dx_{li}  (\frac{\partial T_l}{\partial x_l}) + e_{li}, \quad i = 1, \cdots, n, \label{gastaylor}
\end{eqnarray}
where $e_{gi}$ and $e_{li}$ are the errors in the Taylor's expansion at 
the interface point $x_I$,  $dx_{li} = x_{li} - x_I$, and  
$dx_{gi} = x_{gi} - x_I$.  The continuity of the temperature field $(T_l = T_g = T_I)$  is reflected in the 
first terms of the right-hand side of (\ref{liqtaylor}) and (\ref{gastaylor}). 
In addition to that, we add the 
jump condition (\ref{compnshfluxjump}) in the system of equations 
(\ref{liqtaylor}) and (\ref{gastaylor}) as a constraint. 
Then we obtain the interface temperature $T_I$ 
after solving the  $m + m + 1$ system of equations from minimizing 
the errors \cite{IT02, TKH09}.
The system of equations can be 
re-expressed in the following matrix form 
\begin{eqnarray}
\left( \begin{array}{ccc}
1 & dx_{g1} & 0 \\
\vdots & \vdots & \vdots \\
1 & dx_{gm} & 0 \\
1 & 0 & dx_{l1} \\
\vdots & \vdots & \vdots \\
1 & 0 & dx_{ln}  \\
0 & \kappa_l & -\kappa_g
\end{array} \right ) 
\left( \begin{array}{c}
T_I \\
\\
\frac{\partial T_g}{\partial x_g} \\
\\
\frac{\partial T_l}{\partial x_l} 
\end{array} \right) 
= 
\left( \begin{array}{c}
T_{g1} \\
\vdots \\
T_{gm} \\
T_{l1} \\
\vdots \\
T_{ln}^l \\
0
\end{array} \right) 
+
\left( \begin{array}{c}
e_{g1} \\
\vdots \\
e_{gm} \\
e_{l1} \\
\vdots \\
e_{ln} \\
e_I
\end{array} \right) 
.
\label{interfacesystem}
\end{eqnarray} 
Here the $e_I$ denotes the error in the continuity of of the heat flux at 
the interface. 
This constraint interpolation gives accurate solutions of the 
diffusion equation with discontinuous diffusion coefficients \cite{IT02}. 
After determining the temperature at the liquid interface, we use them 
as Dirichlet boundary conditions for 
the thermodynamic equation (\ref{heateqn}) and obtain the temperature 
for the liquid phase. 

This provides the coupling conditions from the gas phase into the 
liquid phase. To summarize the above descriptions, 
we present the following algorithm for coupling the 
compressible and the incompressible Navier-Stokes equations. 

\subsubsection{Coupling Algorithm I}
 \noindent (i) Generate initial particles with flags as liquid and gas and prescribe 
 the initial values. \\
(ii) Determine the interface particles $x_L$ and $x_R$ from the set of liquid particles. \\
(iii) Solve the compressible Navier-Stokes equations with the 
Dirichlet boundary conditions at the interface (liquid) particles. \\
(iv) Interpolate the pressure and the temperature to the interface 
particles. \\
(v) Solve the incompressible Navier-Stokes equations. \\
(vi) Move all particles with their velocities. \\
(vii) Add or remove particles if necessary, see subsection \ref{fpm}. \\
(viii) Goto (ii) and repeat until the final time is reached. 

\subsection{Coupling of the Boltzmann and the incompressible Navier-Stokes equations}
\label{BoltzNScouple}

The coupling procedure for the Boltzmann and the incompressible 
Navier-Stokes equations is different from the previous one. The Boltzmann 
equation is a mesh-based method since gas molecules have to be sorted 
in cells for the intermolecular collisions. Moreover, we  
sample and store the macroscopic quantities at the cell centers. 
As already described, the incompressible Navier-Stokes equations 
are solved by a meshfree method.  
So, we need to couple the mesh-based and the meshfree method. 
In this coupling procedure, we first generate regular cells in the entire 
domain. In addition to that, we fill the liquid domain by liquid particles. 
The liquid particles overlap with the Boltzmann cells. 
We generate gas molecules outside of the liquid domain. Then we initialize 
the velocity distribution as a Maxwellian 
with the initial parameters $\rho_g(0,x), {\bf u}_g(0,x), T_g(0,x)$ in 
all Boltzmann cells. 
As a consequence we 
have cells with and without gas molecules. If a cell contains no gas molecule 
we deactivate it, while cells with gas molecules are active cells. 
We consider a fixed number of cells for the Boltzmann 
solver. If the cell size is larger than the mean free path, we perform 
a refinement such that the new cell size is smaller than the 
mean free path, see  \cite{TKH09} for details.  

Moreover, we prescribe 
initial conditions for the liquid particles. As in the earlier case 
we can determine the left and right boundaries of the liquid domain. 
First, we solve the Boltzmann equation. After the free flow step in the 
Boltzmann solver we apply the boundary conditions. If gas molecules 
cross the interface to the liquid, they lose the memory of their velocity 
and are reflected diffusively with a Maxwellian distribution having the local 
velocity and temperature parameters of the interface particles. 
This gives the coupling conditions from liquid phase into the gas phase. 

Similarly, if the domain boundaries at $a$ and $b$ are considered as solid walls and 
gas molecules cross these boundaries, they also lose 
the memory of their velocity and reflect diffusively with a 
Maxwellian distribution having the wall temperature and velocity parameters. 

The coupling procedure from the gas phase into the liquid phase is similar to the 
coupling of the compressible and incompressible Navier-Stokes equations. 
At the end of a 
Boltzmann solver step we compute the stress tensor and heat flux from 
the gas molecules and store them 
at the cell centers. Moreover, we also compute the density, velocity,  
temperature and pressure and store them at the cell centers. Then we approximate the 
pressure at the interface from the 
neighboring activated cell centers of the Boltzmann solver using equation (\ref{boltzpressjump}). 
The temperature at the interface is approximated similarly as in  
subsection \ref{NSNScouple}, with the 
conditions (\ref{tempjump}) and (\ref{boltzhfluxjump}). Note 
that the system of equations differs slightly from (\ref{interfacesystem}) 
in this case and is re-expressed as 
\begin{eqnarray}
\left( \begin{array}{ccc}
1 & dx_{g1} & 0 \\
\vdots & \vdots & \vdots \\
1 & dx_{gm} & 0 \\
1 & 0 & dx_{l1} \\
\vdots & \vdots & \vdots \\
1 & 0 & dx_{ln}  \\
0 & -\kappa_l & 0
\end{array} \right ) 
\left( \begin{array}{c}
T_I \\
\\
\frac{\partial T_g}{\partial x_g} \\
\\
\frac{\partial T_l}{\partial x_l} 
\end{array} \right) 
= 
\left( \begin{array}{c}
T_{g1} \\
\vdots \\
T_{gm} \\
T_{l1} \\
\vdots \\
T_{ln} \\
q_g
\end{array} \right) 
+
\left( \begin{array}{c}
e_{g1} \\
\vdots \\
e_{gm} \\
e_{l1} \\
\vdots \\
e_{ln} \\
e_I
\end{array} \right) 
.
\label{interfacesystem1}
\end{eqnarray}

It is well known that  in all DSMC type solvers  
there are some statistical fluctuations in the solutions. If we 
pass these fluctuating data, they destabilize the Navier-Stokes solver. 
Therefore, we need a smoothing operator. Here we have used the Shepard interpolation.
For example, for a function $\psi$ at a cell center $x$, 
the Shepard interpolation 
is defined as 
\begin{equation}
\tilde \psi(x) = \frac{\sum_{i=1}^m w_i \psi(x_i)}{\sum_{i=1}^m w_i},
\label{shephard}
\end{equation}
where $m$ is the number of neighboring  active cell centers $x_i$, 
and $w_i$ are the weight functions depending on the 
distance between $x$ and $x_i$. We have chosen a Gaussian weight function 
given by 
\begin{eqnarray} 
w_i = w( x_i - x; h) =
\left\{ 
\begin{array}{l}  
exp (- \alpha \frac{( x_i - x  )^2 }{h^2} ), 
\quad \mbox{if    }  \frac{| x_i - x  |}{h} \le 1 
\\
 0,  \qquad \qquad \quad \quad \quad \mbox{else}
\end{array}
\right.
\label{weight}
\end{eqnarray}
with $h$ being the radius of interaction which is three times the 
cell distance. 
$ \alpha $ is a  positive constant, taken as $2$, which, however, can be varied. 
We have already used this type of smoothing for coupling the Boltzmann and 
the Navier-Stokes equations \cite{TKH09}, with the result that for small 
Knudsen numbers the locally smoothed Boltzmann solutions 
are compatible with the non-smoothed Navier-Stokes solutions. 
This means that there is negligible smearing of the solutions with this 
type of local smoothing.  

After approximating the pressure and temperature at the gas-liquid interface, we 
solve the incompressible Navier-Stokes equations in a similar manner 
as in subsection \ref{NSNScouple}. This completes the description of the
coupling from the Boltzmann domain into the incompressible 
Navier-Stokes domain.

Recall that in the liquid phase a Lagrangian particle method is applied. 
Therefore, 
after moving the liquid particles, we may have some gas molecules inside 
the liquid domain. Since we switch to the Boltzmann solver after a Navier-Stokes step has been completed, 
we do not reflect them from liquid interface at this 
level, but we treat them after the free-flow step of the Boltzmann solver.  

Summarizing the above we present the following algorithm for the coupling of 
the Boltzmann and the incompressible Navier-Stokes equations. 

\subsubsection {Coupling Algorithm II} 
 \noindent (i) Generate Boltzmann cells in the entire domain and generate 
liquid particles which overlap the Boltzmann cells. \\
(ii) Generate gas molecules outside the 
liquid domain according to a Maxwellian distribution with the initial 
parameters and prescribe the initial conditions for the liquid particles. \\
(iii) Determine the interface particles $x_L$ and $x_R$ for the liquid phase. \\
(iv) Solve the Boltzmann equation with the gas-liquid interface taking the role of a 
moving wall.\\
(v) Compute the moments in the Boltzmann cells and 
smooth them according to (\ref{shephard}).\\
(vi) Interpolate the pressure and temperature to the interface 
particles. \\
(vii) Solve the incompressible Navier-Stokes equations. \\
(viii) Move the liquid particles with their velocities. \\
(ix) Add or remove liquid particles, if necessary, see subsection \ref{fpm}.\\
(x) Goto (iii) and repeat until the final time is reached.

\section{Numerical Tests}
\label{numsec}
As we have already mentioned, we consider the computational domain  
$\Omega = [a,b]$, where $a = 0$ in all cases 
and $b$ varies. In $\Omega_l$ we solve the incompressible Navier-Stokes 
equations whereas in $\Omega_g$ we solve either the Boltzmann equation or the 
compressible Navier-Stokes equations. 
Then we compare the results obtained from both coupling algorithms. 
We consider three different test cases. In all cases a 
liquid droplet is surrounded by a compressible gas inside the computational 
domain, where the boundary points $a$ and $b$ always belong 
to $\Omega_g$. While coupling the compressible and incompressible 
Navier-Stokes equations we consider a total number 
of initial particles $N=200$ for both liquid and gas domains. 
The corresponding initial grid spacing is $\Delta x = 1/N$. To determine 
the neighboring particles in the FPM, the interaction radius $h$ is taken equal to 
$3\times\Delta x$.

For the coupling of the Boltzmann and the 
incompressible Navier-Stokes equations we generate $N = 200$ regular cells in 
$\Omega$ and add liquid particles which cover the liquid 
domain $\Omega_l$. Note that in algorithm I, $\Omega_g$ and $\Omega_l$ are 
disjoint, but same is not true in algorithm II. 
The initial grid spacing of liquid particles is again equal to $1/N$.  
The Boltzmann cell size is refined according to the size of the 
mean free path, see \cite{TKH09} for details. 
As we have mentioned initially, the 
hard-sphere model for the collision cross section is employed for the 
Boltzmann equation. Moreover, 
we have used the following parameters for all test cases. 
The gas is Argon with a molecular mass $m = 6.63\times 10^{-26}~kg$. 
For the Boltzmann constant we have $k=1.38\times 10^{-23}~J/K$, for the 
molecular diameter $d=3.68\times 10^{-10}~m$, and we obtain 
the gas constant $R = 208~J/(kg K)$. The dynamic viscosity and 
thermal conductivity for the compressible Navier-Stokes equations 
are assumed to be constant and are evaluated with the initial temperature according to (\ref{mukappa}). 

For the liquid phase, in all three cases we assume an 
initial temperature equal to $298~K$. The thermal conductivity 
and the specific heat capacity of the liquid are taken to be those of 
water, giving values of $\kappa_l = 0.63~J/(mKs)$ and 
$c_p = 4181~J/(kgK)$, respectively. 

\subsection{ Test 1 }
We first study a liquid droplet with nonzero initial velocity. 
This test case has been considered by Caiden et al \cite{CFA00} for the 
inviscid case. We first consider the interval 
$\Omega = [0~m,1\times 10^{-4}~m]$. Initially, a liquid droplet occupies the 
domain $\Omega_l = [4\times 10^{-5}~m, 6\times 10^{-5}~m]$, 
while the gas occupies the rest of the domain. 
The gas has the initial values $\rho_g = 1.226~kg/m^3, u_g = 0~m/s$ and 
$p_g = 1\times10^5~Pa$. The temperature is obtained from the 
equation of state. 
The initial values for the liquid phase are $\rho_l = 1000~kg/m^3, 
u_l = 100~m/s, T_l = 298~K$. The initial pressure in the liquid is taken to be the 
same as the initial pressure in the gas, however, this will change due to the coupling of the phases. 

The boundary values for the compressible Navier-Stokes equations are given by 
$u_g(t,a) = u_g(t,b) = u_g(0,x)$ and $T_g(t,a) = T_g(t,b) = T_g(0,x)$. 
Moreover, for the Boltzmann solver, if gas molecules cross the physical 
domain, we reflect them back in the same was as in the case of 
the gas-liquid interface boundary conditions, using the boundary velocity and temperature 
$u_g(t,a) = u_g(t,b) = u_g(0,x)$ and $T_g(t,a) = T_g(t,b) = T_g(0,x)$. 

The characteristic length is taken to be the length of the droplet which is 
equal to $2\times 10^{-5}~m$. The corresponding Knudsen 
number is $\epsilon = 0.0045$.   
The simulations are stopped after $t=5.2\times 10^{-8}~s$. 
\begin{figure}
\begin{center}
\epsfig{file=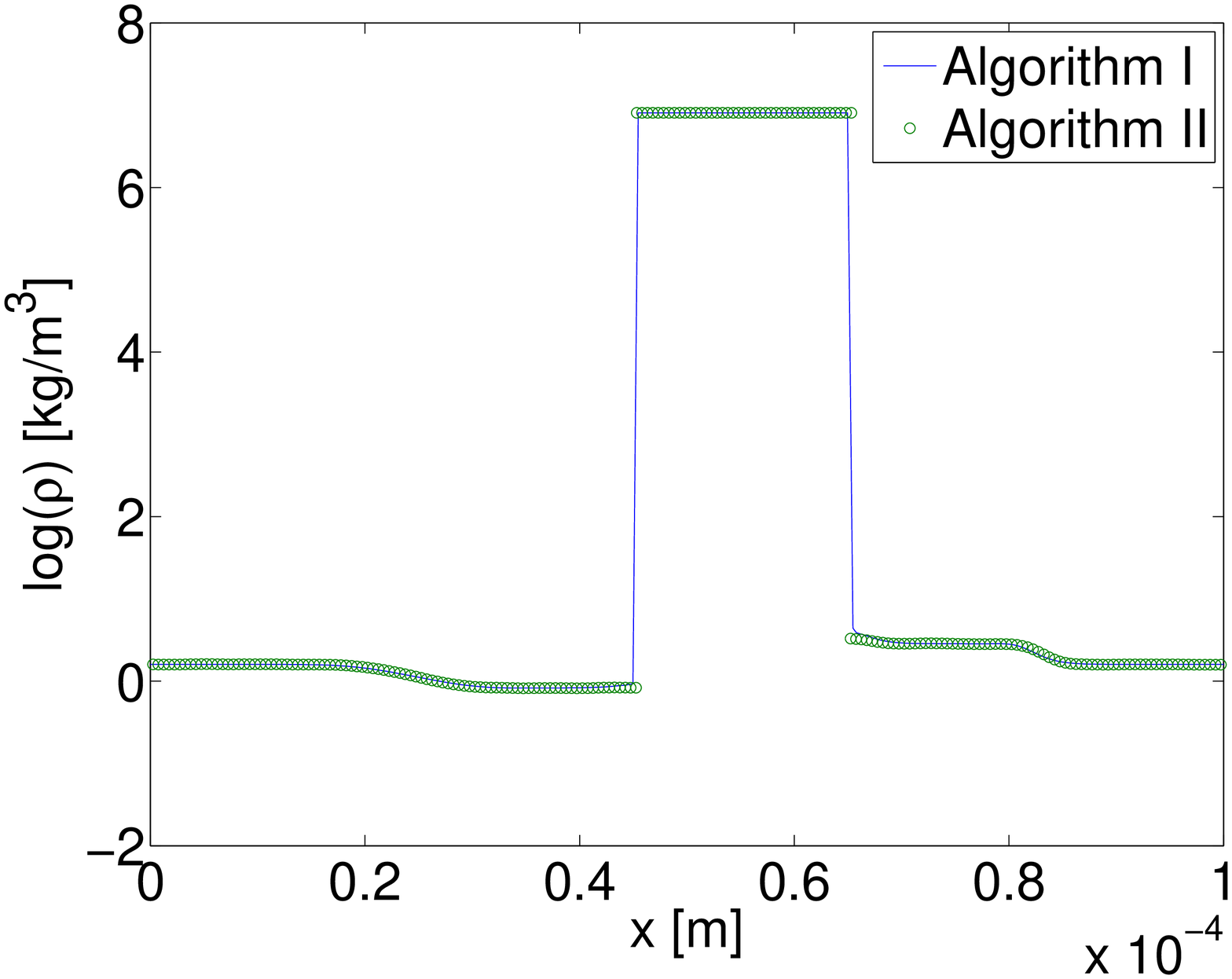,width=6.85cm}\hfill
\epsfig{file=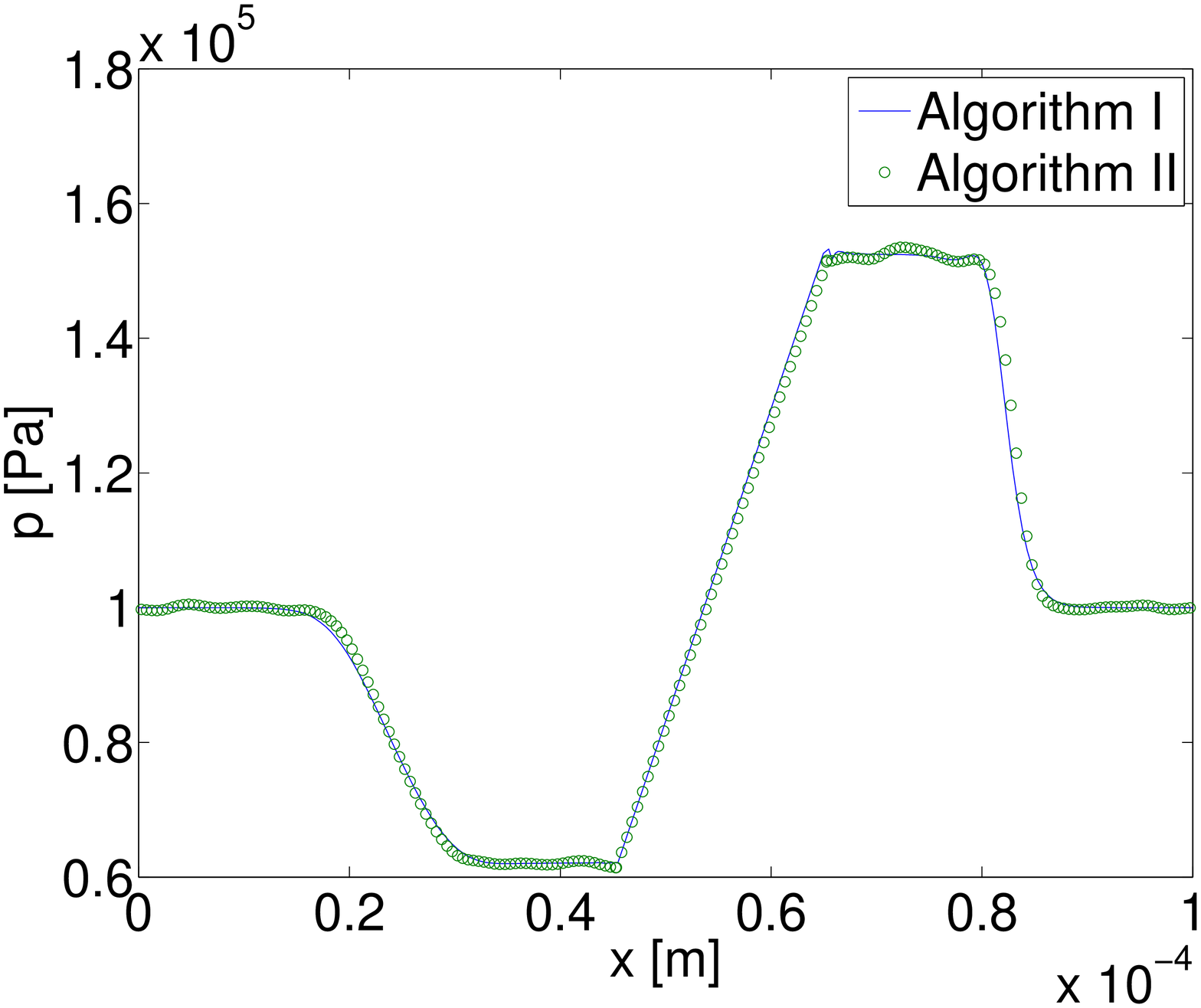,width=6.85cm}
\epsfig{file=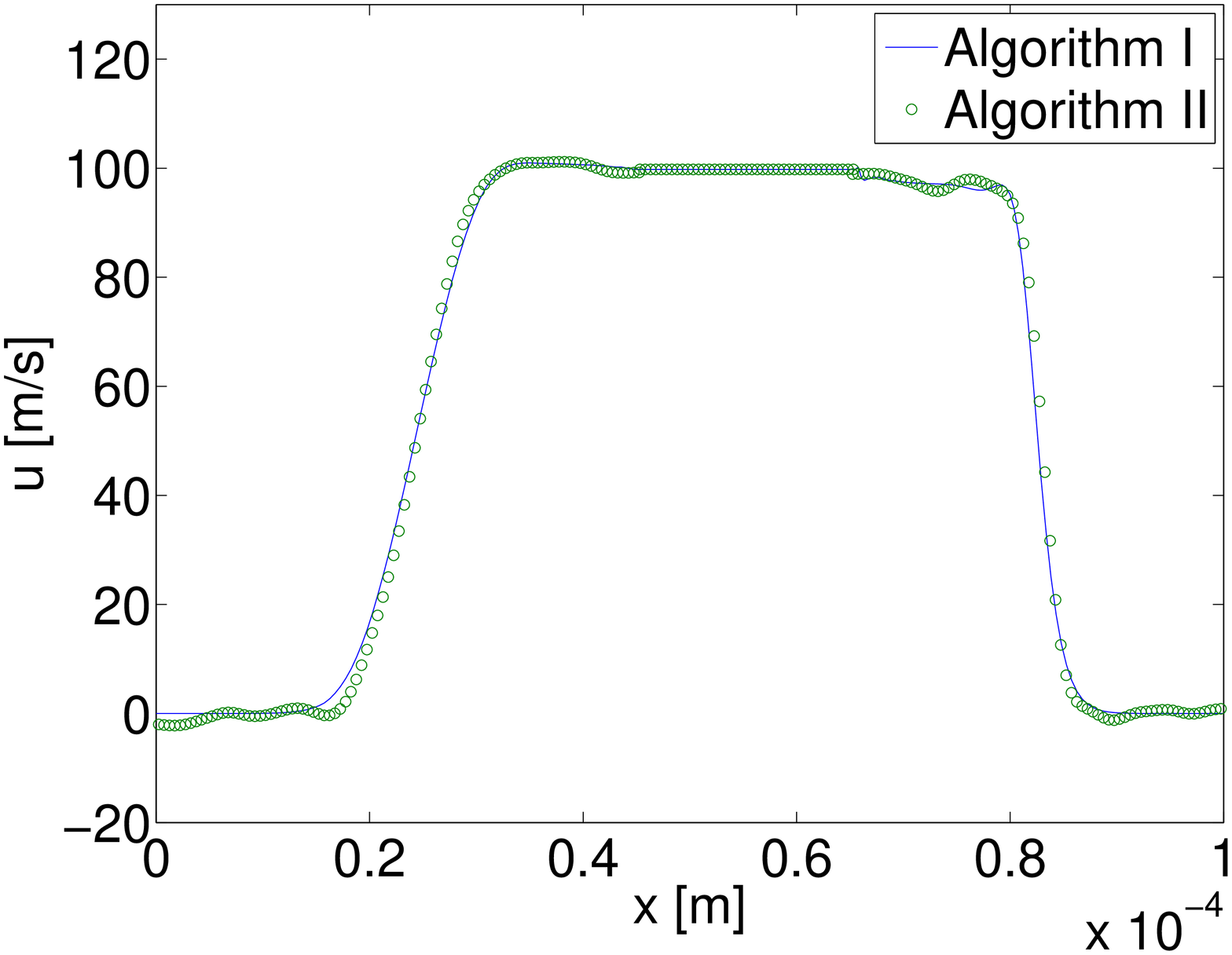,width=6.85cm}\hfill
\epsfig{file=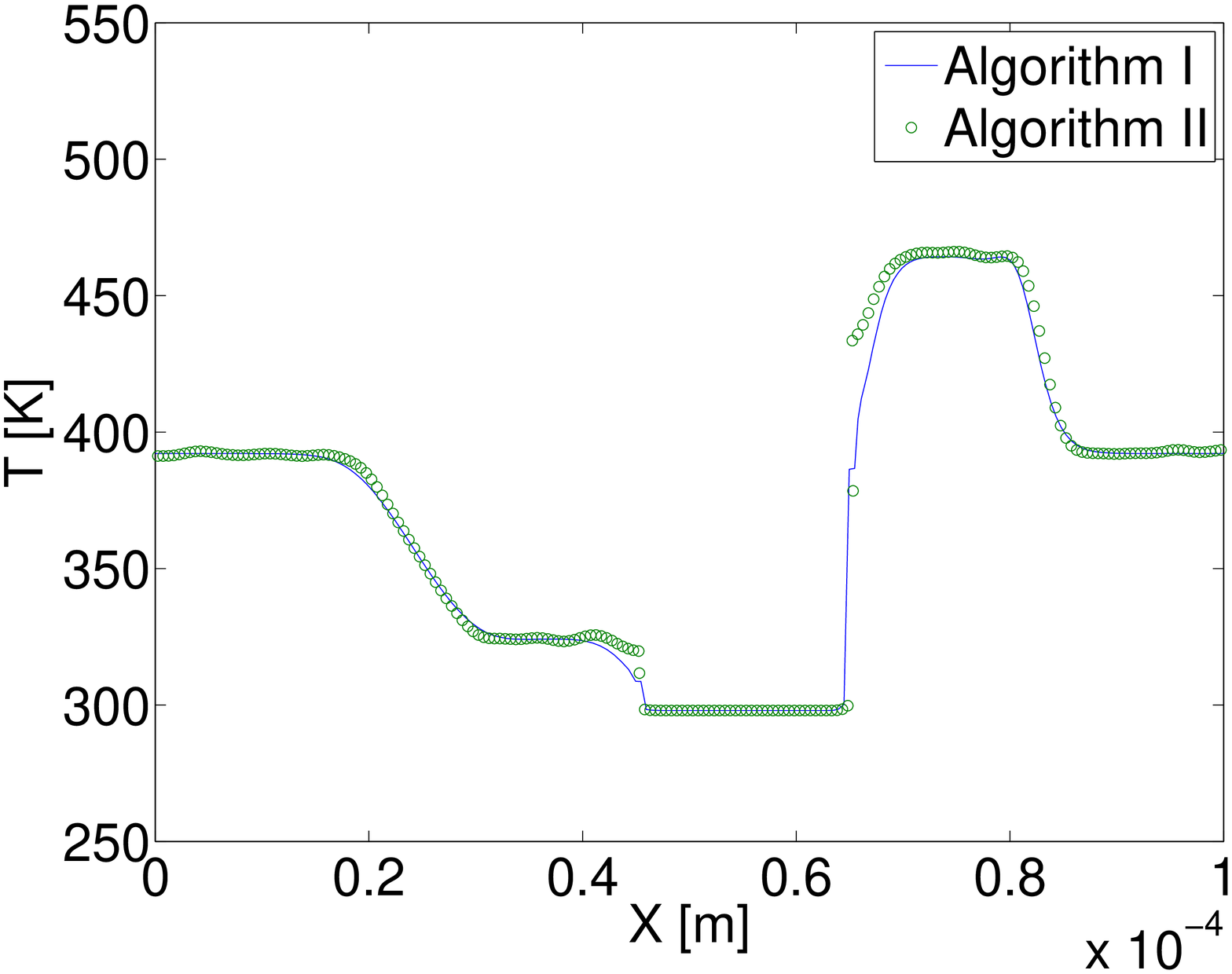,width=6.85cm}
 \end{center}
\caption{  \em Test1: Logarithm of density(top left), pressure(top right), velocity(bottom left) and temperature(bottom right) for a Knudsen number of 0.0045. The solid lines are results from Algorithm I, the circles are from Algorithm II. }
\label{test1case1}
\end{figure}

The initial number of gas molecules is $5000$ per 
cell of size $\Delta x$ for the simulation of the Boltzmann equation. 
We observe in figure \ref{test1case1} that the solutions obtained from the 
coupling of the Boltzmann and the incompressible Navier-Stokes equations (Algorithm II) 
match with the solutions obtained from the coupling of the compressible 
and incompressible Navier-Stokes equations (Algorithm I). 
 
As in \cite{CFA00} we notice that 
the liquid droplet moves to the right causing a compression wave in the gas 
ahead of it and an expansion wave behind it. Moreover, the 
gas is heated ahead of droplet and cooled behind it.
To a very good approximation, the pressure profile inside the droplet is linear. This is expected and can be explained by the
following argument. The liquid is decelerated by compressing the gas volume at its right. By the equivalence principle of general relativity, an
acceleration (or deceleration) is equivalent to the action of a gravitational field. Therefore a linear hydrostatic pressure profile as in a quiescent
liquid in a gravitational field is found. 
\begin{figure}
\begin{center}
\epsfig{file=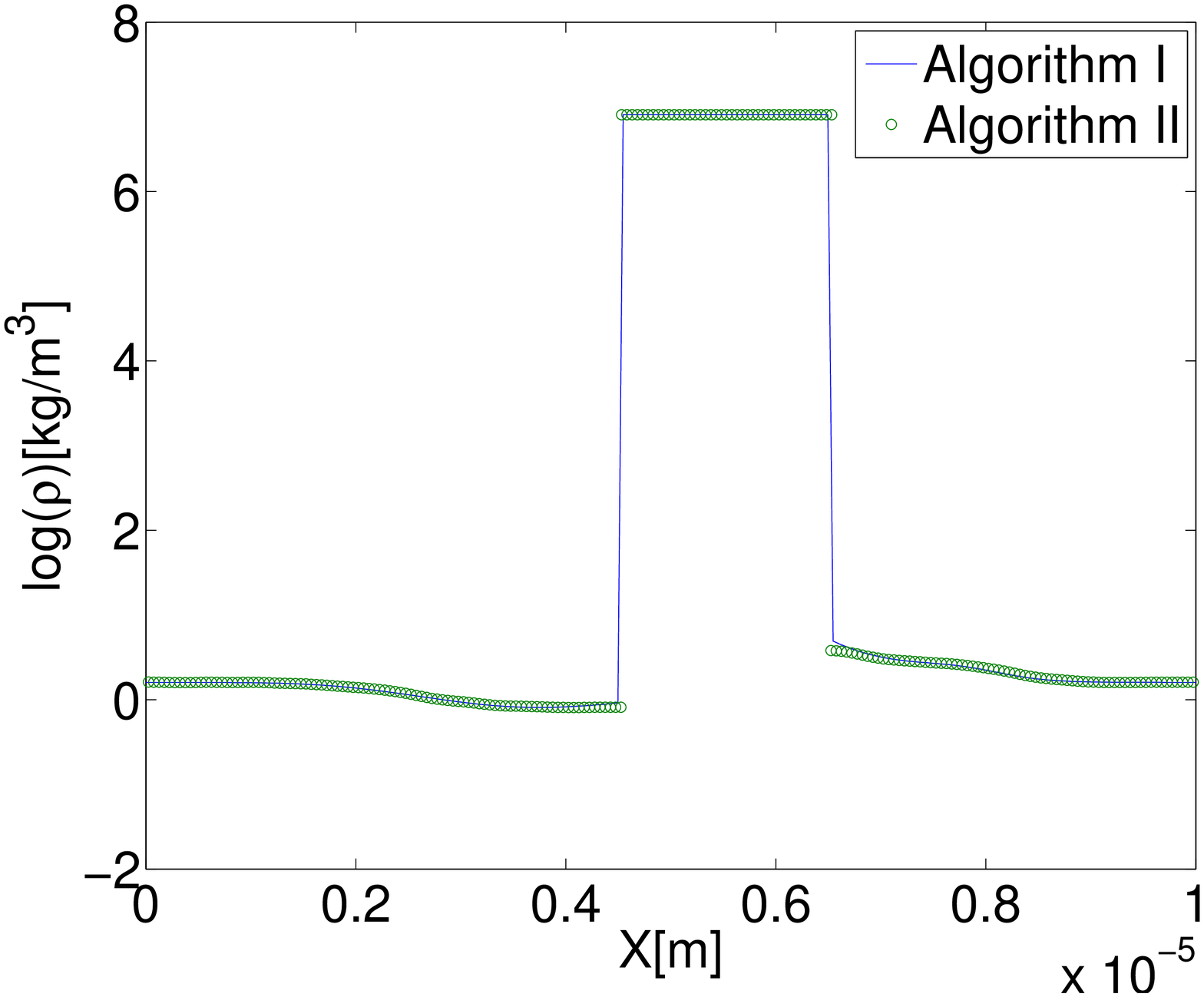,width=6.85cm}\hfill
\epsfig{file=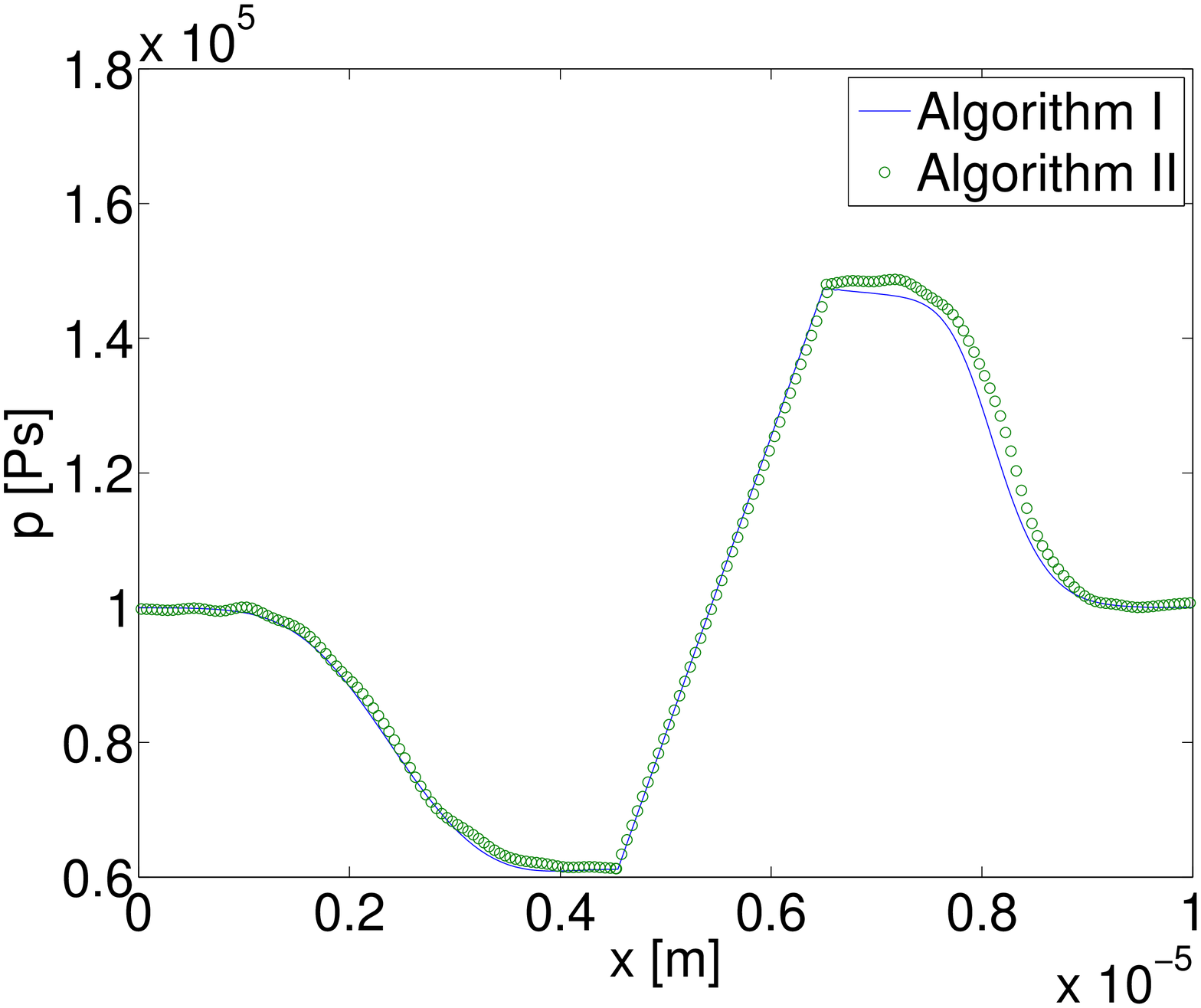,width=6.85cm}
\epsfig{file=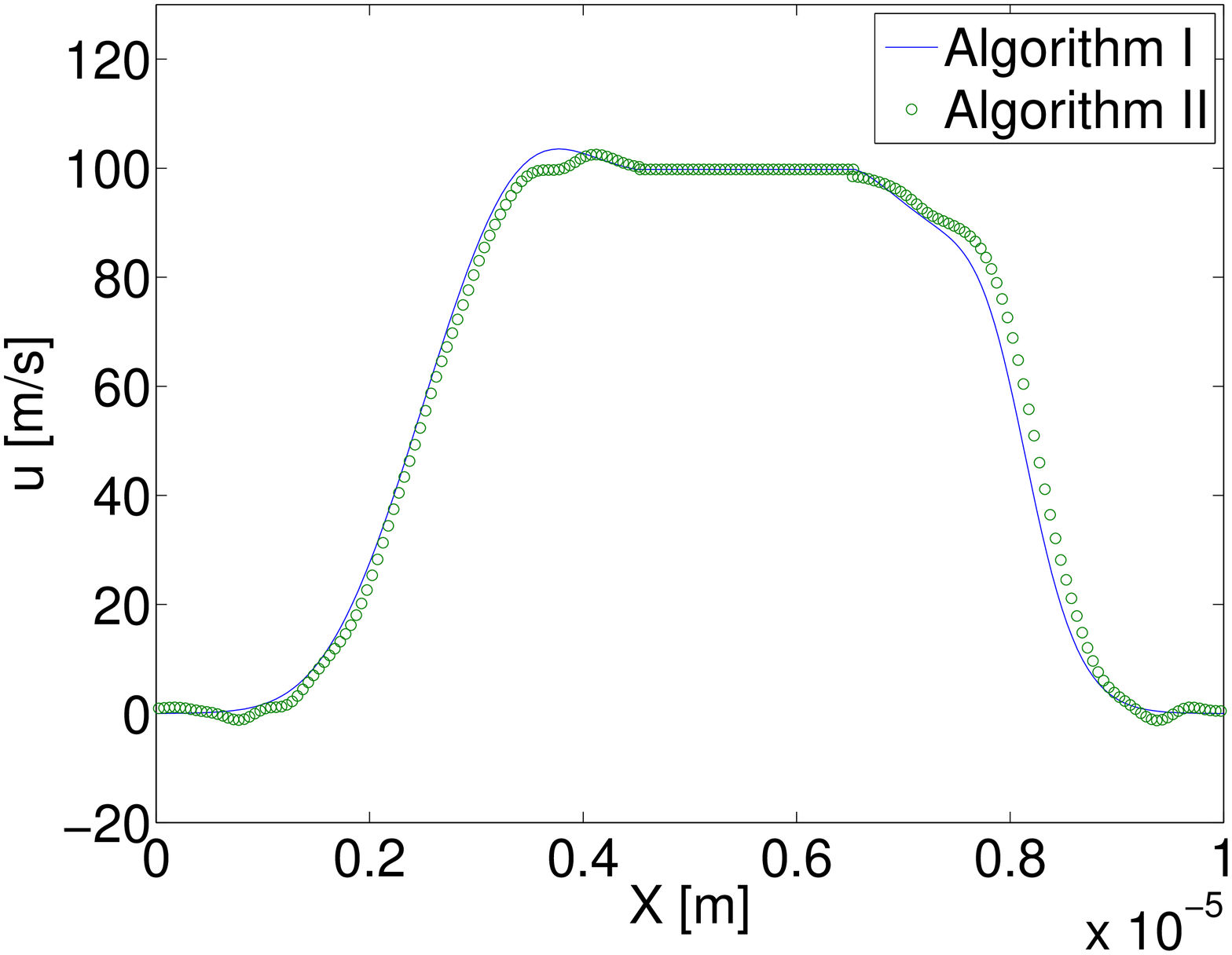,width=6.85cm}\hfill
\epsfig{file=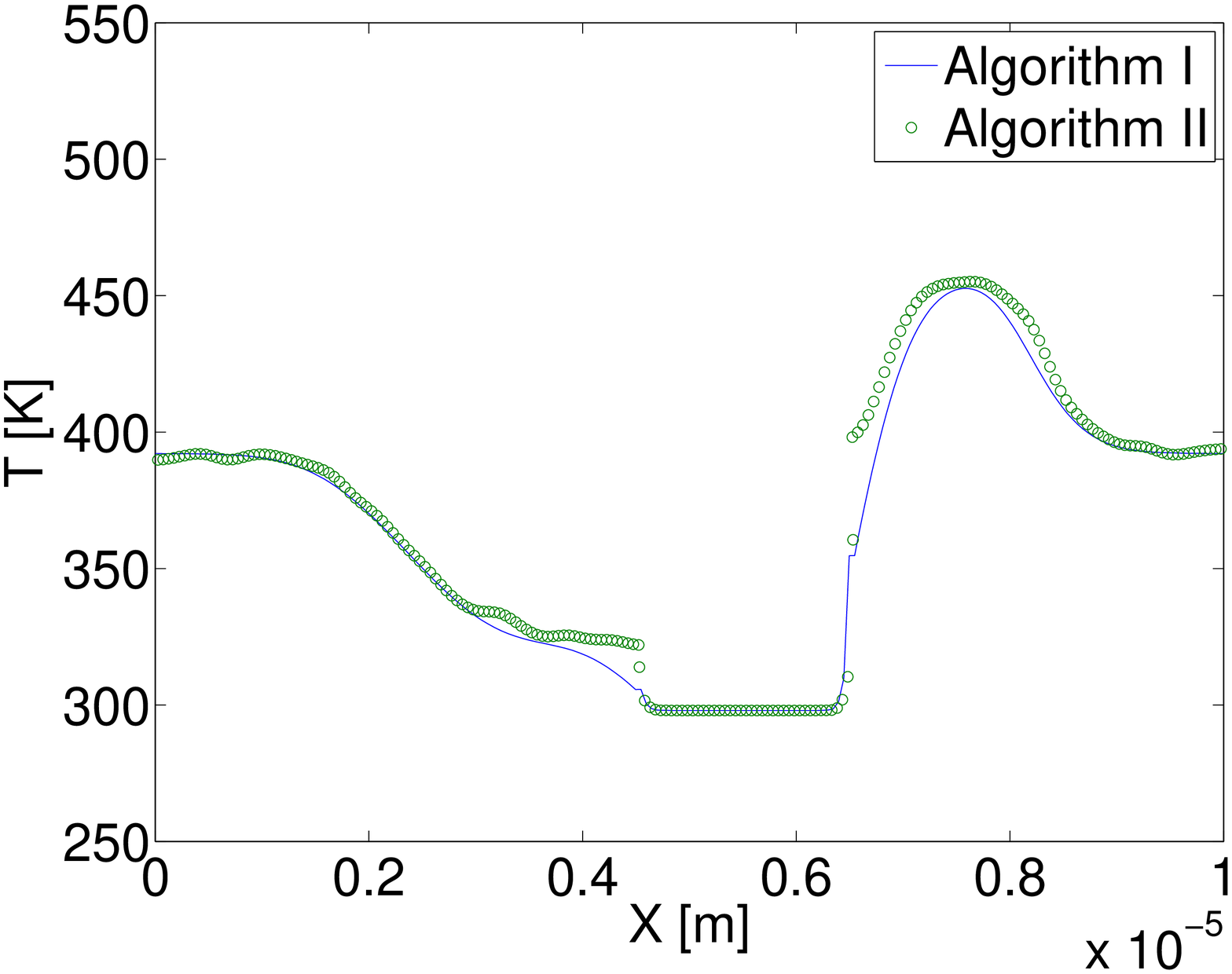,width=6.85cm}
 \end{center}
\caption{ \em Test1: Logarithm of density(top left), pressure(top right), velocity(bottom left) and temperature(bottom right) for a Knudsen number of 0.045. The solid lines are results from Algorithm I, the circles are from Algorithm II. }
\label{test1case2}
\end{figure}

Next, we consider a $10$ times smaller domain 
$\Omega = [0,1\times 10^{-5}~m]$ and the same initial values as above. Initially 
a liquid droplet occupies the domain   
$\Omega_l = [4\times 10^{-6}~m, 6\times 10^{-6}~m]$, while the gas 
occupies the rest of the domain. The characteristic 
length is taken as the size of the droplet, which is equal to 
$2\times 10^{-6}$, giving a 
Knudsen number $\epsilon = 0.045$. 
The simulation is stopped after a time $t = 5.2\times 10^{-9}~s$.
The results are plotted in figure \ref{test1case2} and show similar trends
as the curves in figure \ref{test1case1}. However, at such values of the Knudsen number first deviations from the results of
continuum theory become visible. 

\begin{figure}
\begin{center}
\epsfig{file=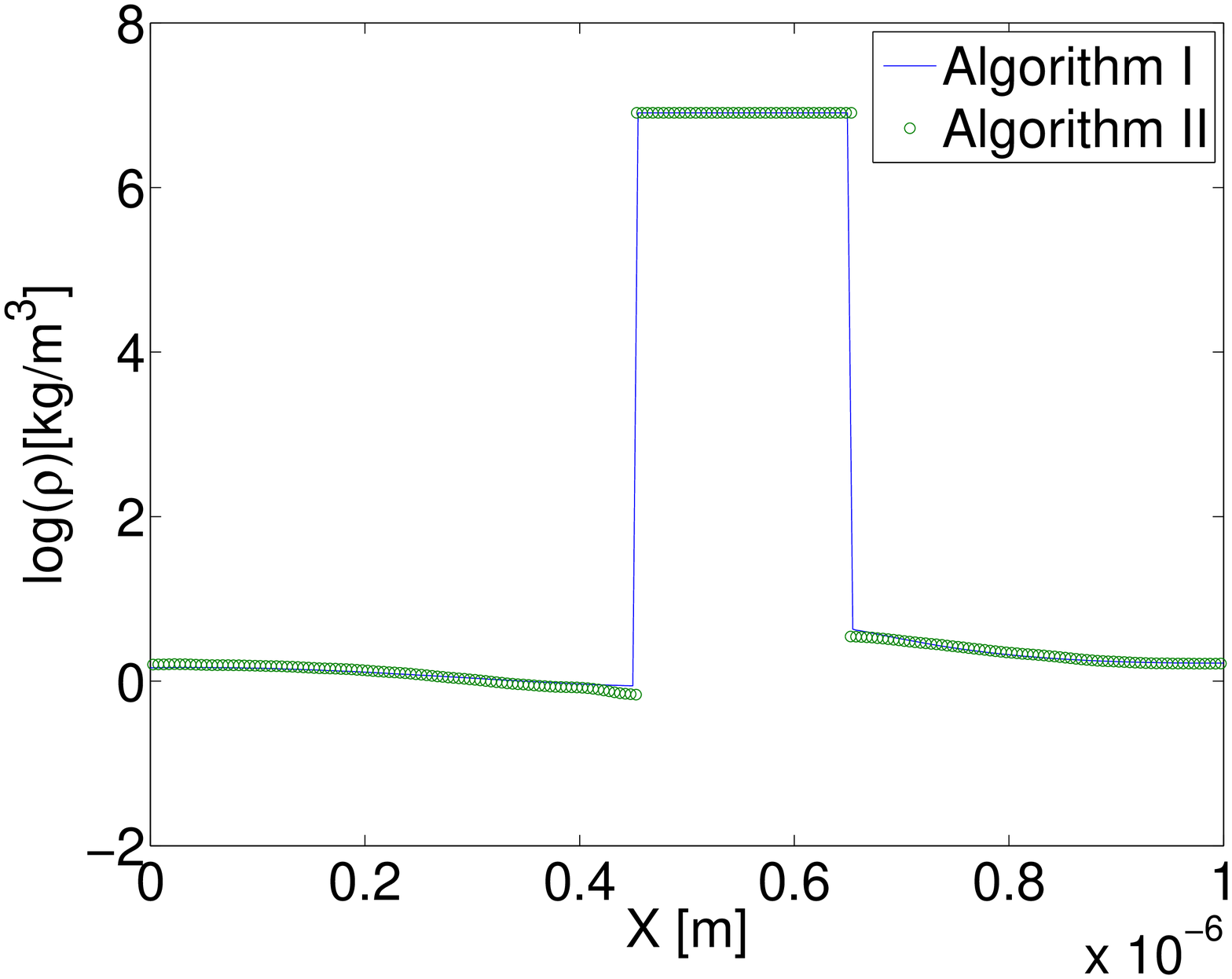,width=6.85cm}\hfill
\epsfig{file=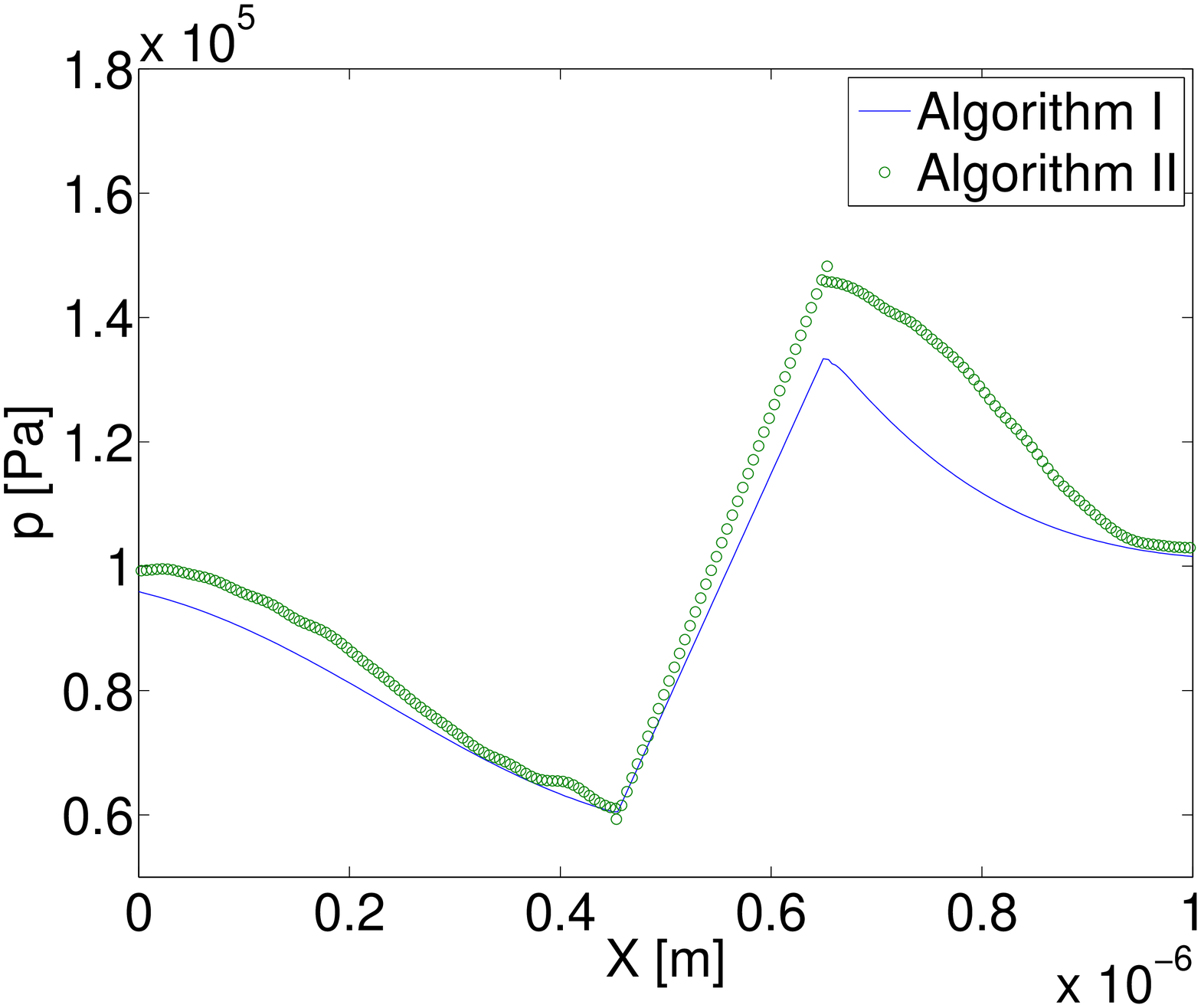,width=6.85cm}
\epsfig{file=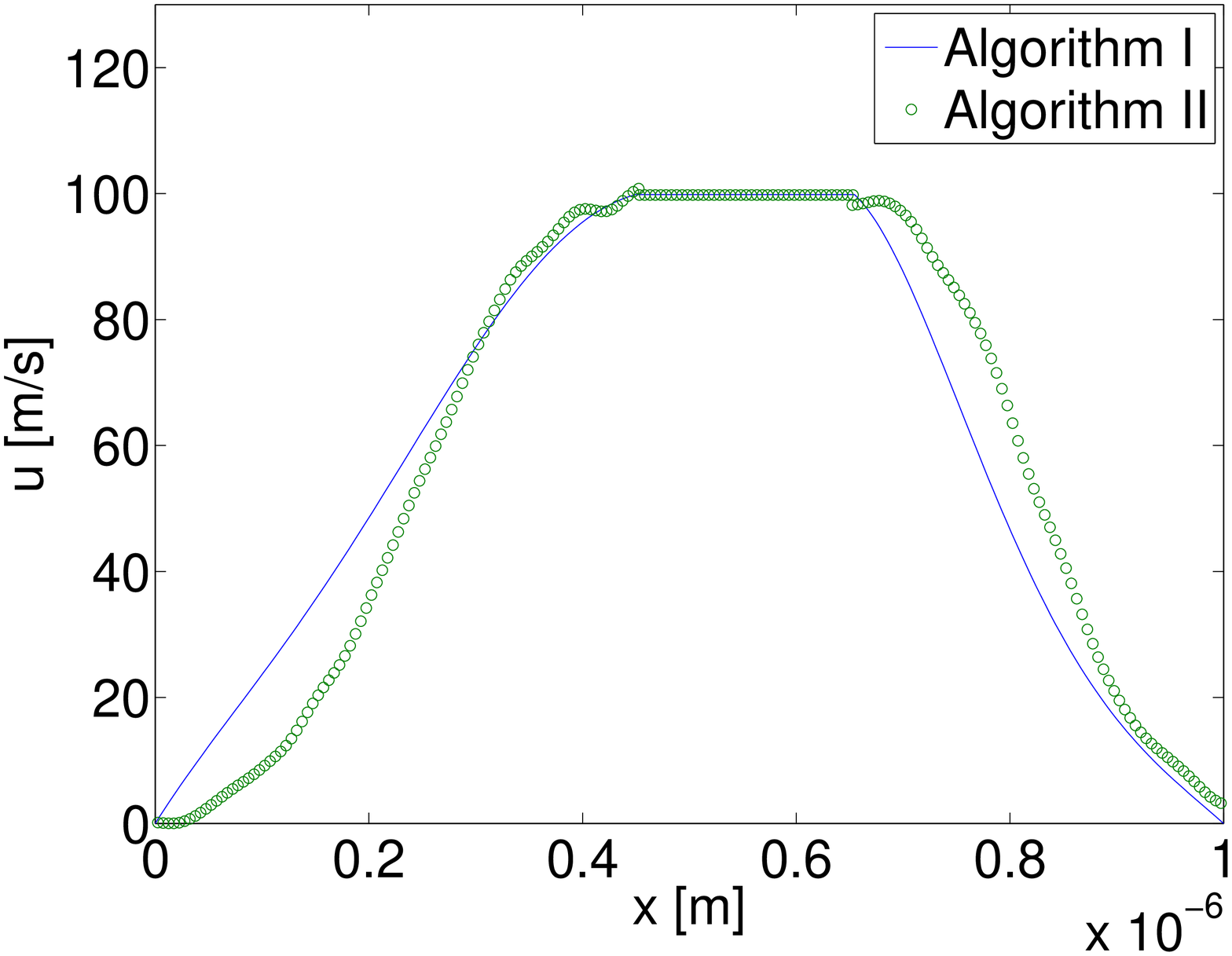,width=6.85cm}\hfill
\epsfig{file=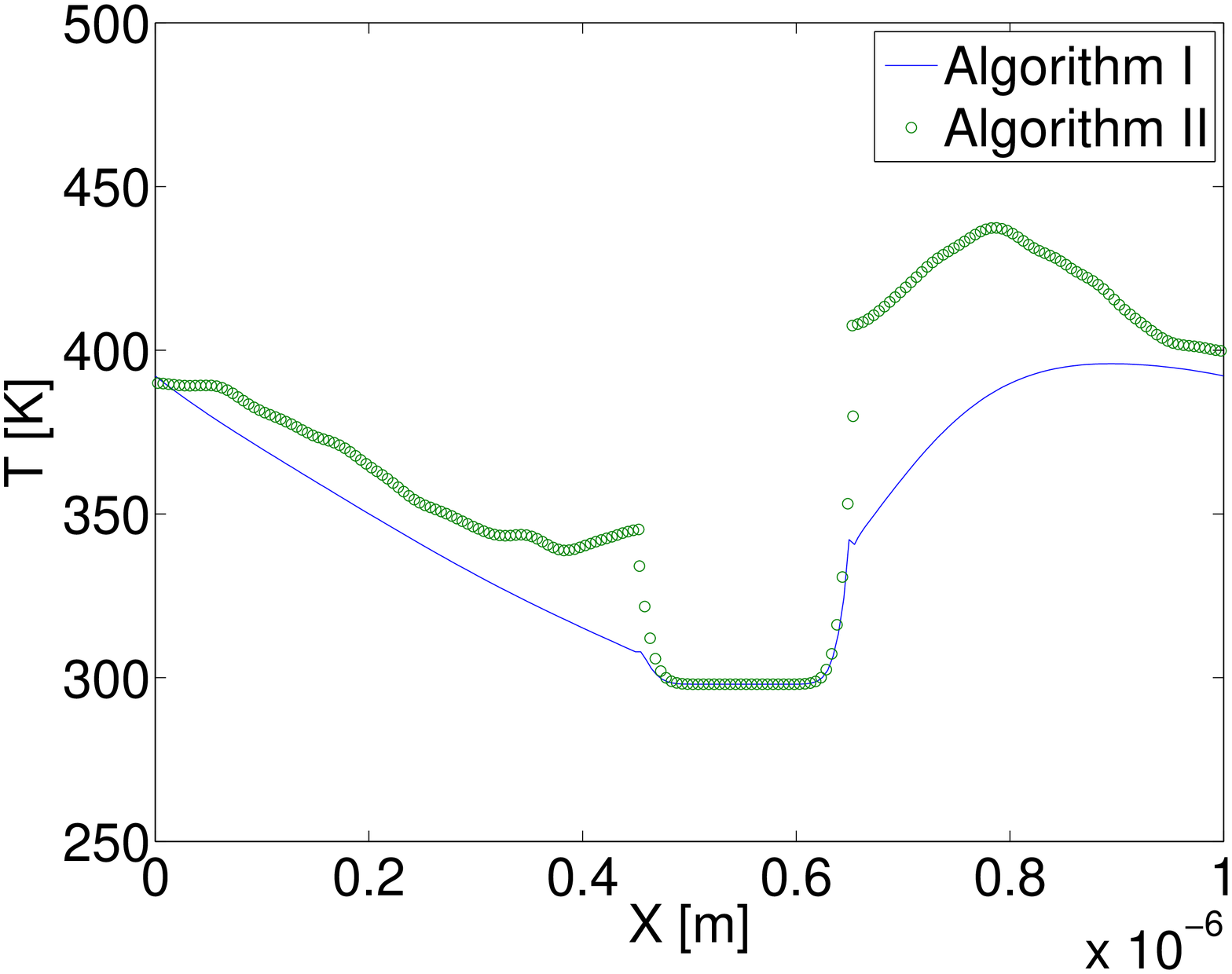,width=6.85cm}
 \end{center}
\caption{ \em Test1: Logarithm of density(top left), pressure(top right), velocity(bottom left) and temperature(bottom right) for a Knudsen number of 0.45. The solid lines are results from Algorithm I, the circles are from Algorithm II. }
\label{test1case3}
\end{figure}

We have further increased the Knudsen number by decreasing the computational domain and characteristic length by a factor of 10. 
In this case we have $\epsilon = 0.45$, and simulation was stopped after 
$t = 5.2\times 10^{-10}~s$. We have plotted the results in 
figure \ref{test1case3}. 
For this Knudsen number we see that the solutions 
from the coupling algorithm I significantly deviate from 
the ones obtained from coupling algorithm II. This is, as expected, due 
to the failure of the compressible Navier-Stokes equations for larger 
Knudsen numbers. A constant velocity and a linear pressure profile are obtained inside the liquid. 
It also becomes apparent that the solution of the 
Boltzmann equation yields much more pronounced jumps in the 
temperature fields than that of the Navier-Stokes equations, 
similar as in the case 
of Sod's 1D shock tube problem \cite{TKH09}.

\subsection{Test 2}

The second test case is also taken from paper of Caiden et al \cite{CFA00}, 
where the authors have considered only the Euler equations without 
the energy equation. We again consider the 
computational domain $\Omega = [0~m,1\times 10^{-4}~m]$. The liquid initially 
occupies the domain $\Omega_l=[4\times 10^{-5}~m, 6\times 10^{-5}~m]$, and 
the rest of the domain is filled with gas. A shock wave is 
initially located at $x=1\times 10^{-5}~m$, with a post shock state 
$p_g(0,x) = 148407.3~Pa, u_g(0,x) = 89.981~m/s$ and 
$\rho_g(0,x) = 2.214~kg/m^3$. 
On the right of 
 $x=1\times 10^{-5}~m$ the initial state is given as 
$p_g(0,x) = 98066.5~Pa, u_g(0,x) = 0~m/s$ and 
$\rho_g(0,x) = 1.58317~kg/m^3$. 
The initial temperature of the gas is again computed from the equation of state. 
The characteristic length is equal to $2\times 10^{-5}~m$. 
The pre- and post-shock Knudsen numbers are 
$0.00259$ and $0.00696$, respectively. 
The initial conditions for the liquid are $p_l = 98066.5~Pa, u_l = 0~m/s, 
T_l = 298~K$. At first, we first consider a liquid density of $\rho_l = 1000~kg/m^3$. 

The boundary conditions for the compressible Navier-Stokes equations are 
$u_g(t,a) = 89.981~m/s, T_g(t,a) = T(0,a)$,  $T_g(t,b) = T_g(0,b)$, 
and we assume a vanishing velocity gradient at $x=b$.

\begin{figure}
\begin{center}
\epsfig{file=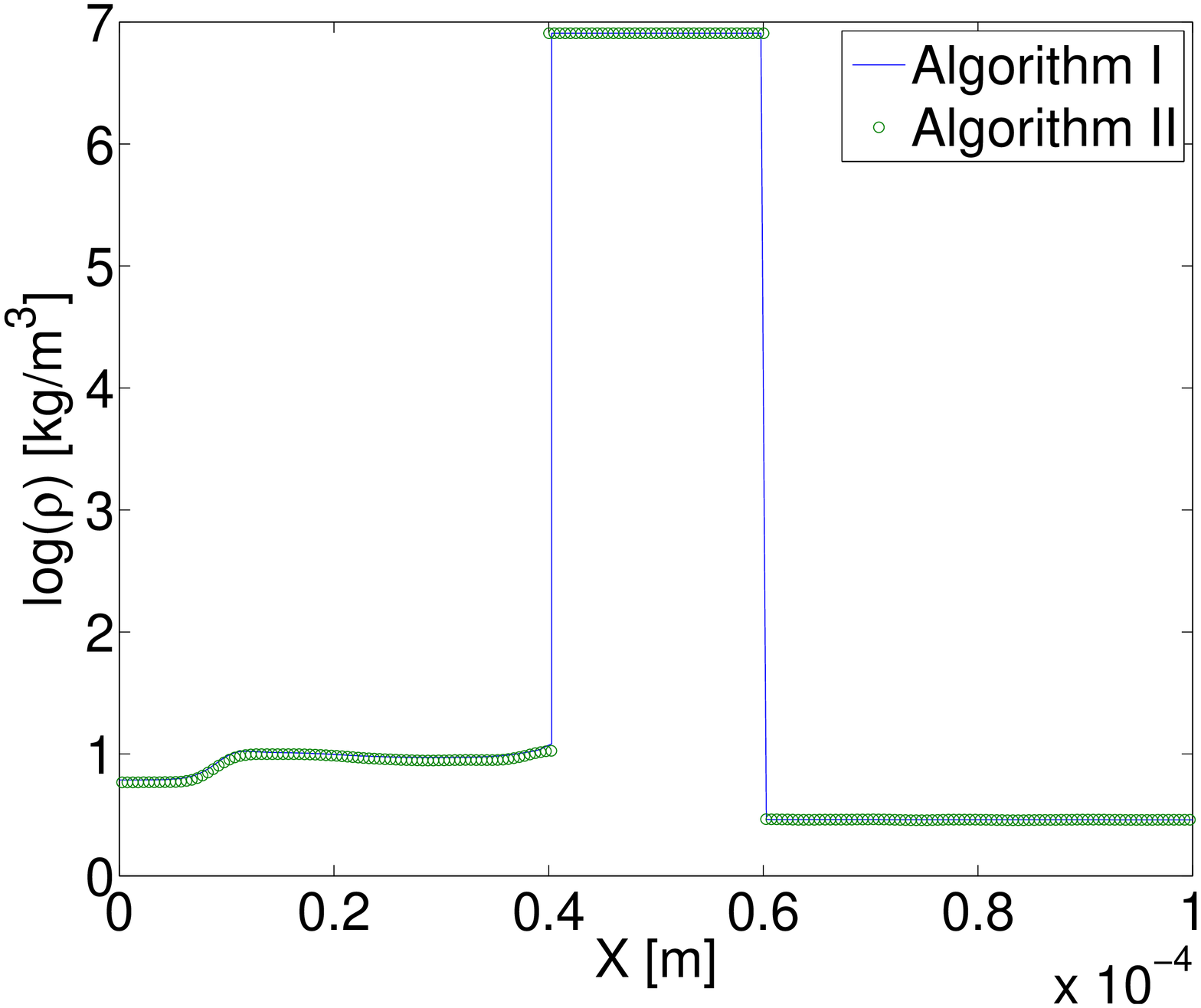,width=6.85cm}\hfill
\epsfig{file=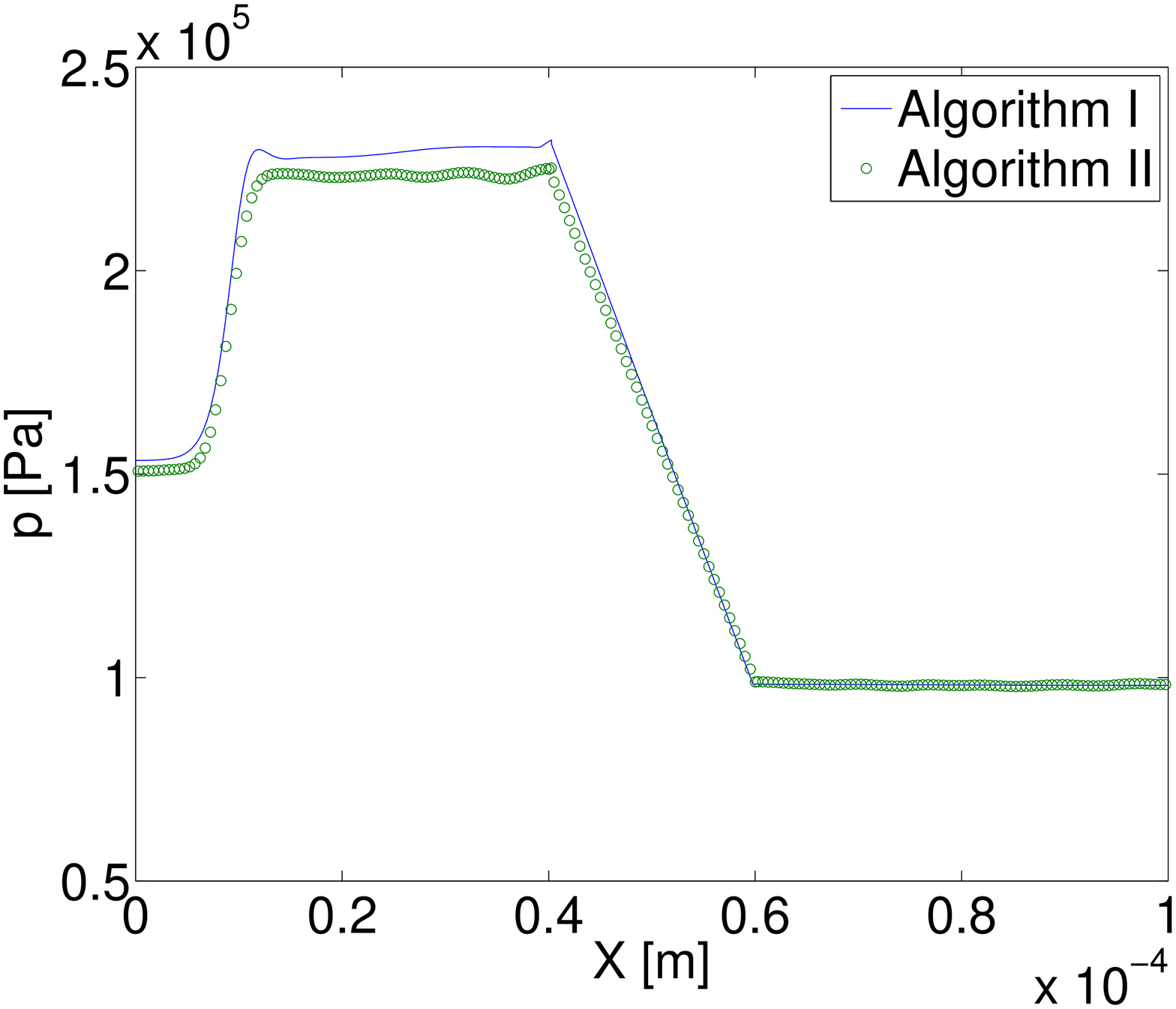,width=6.85cm}
\epsfig{file=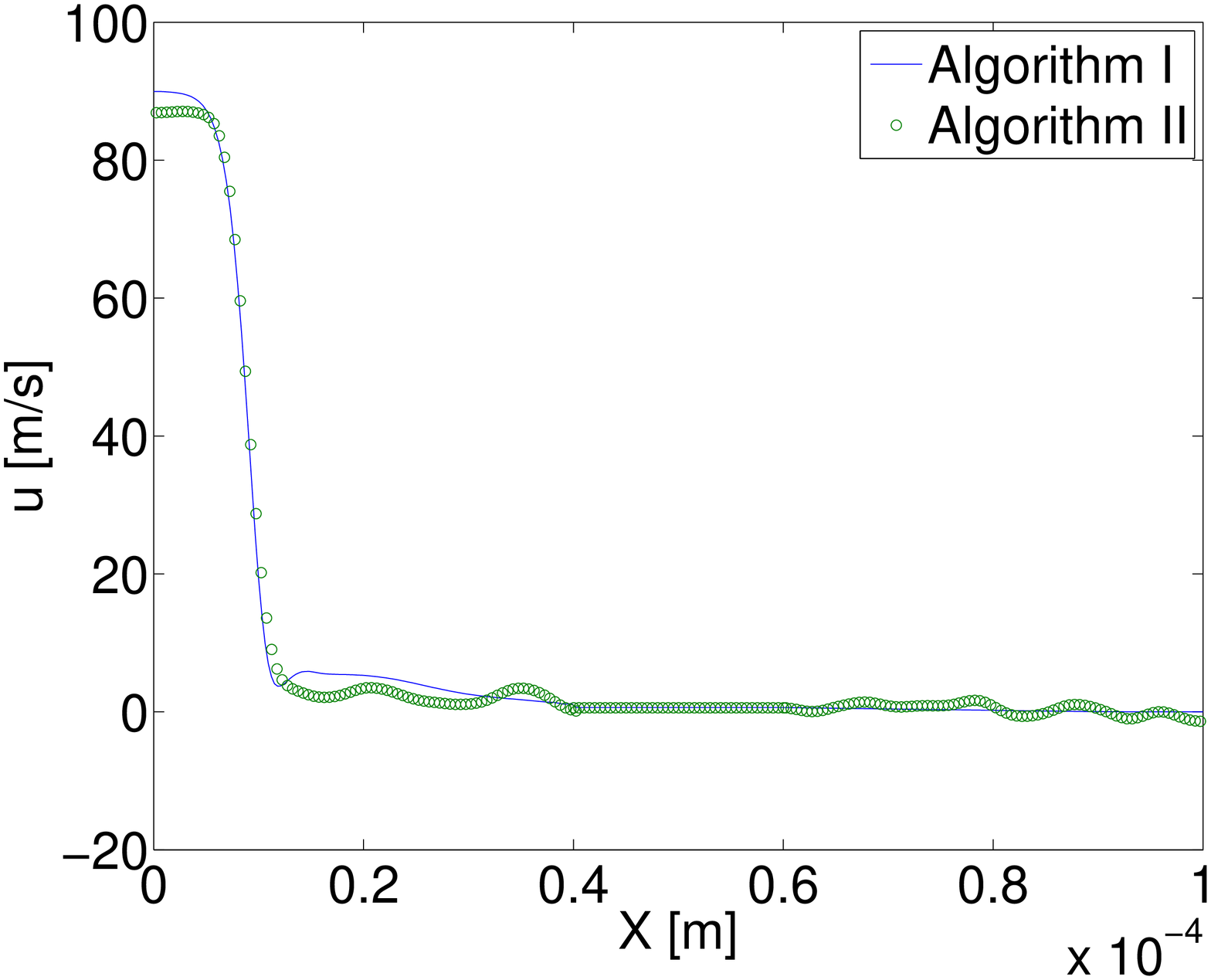,width=6.85cm}\hfill
\epsfig{file=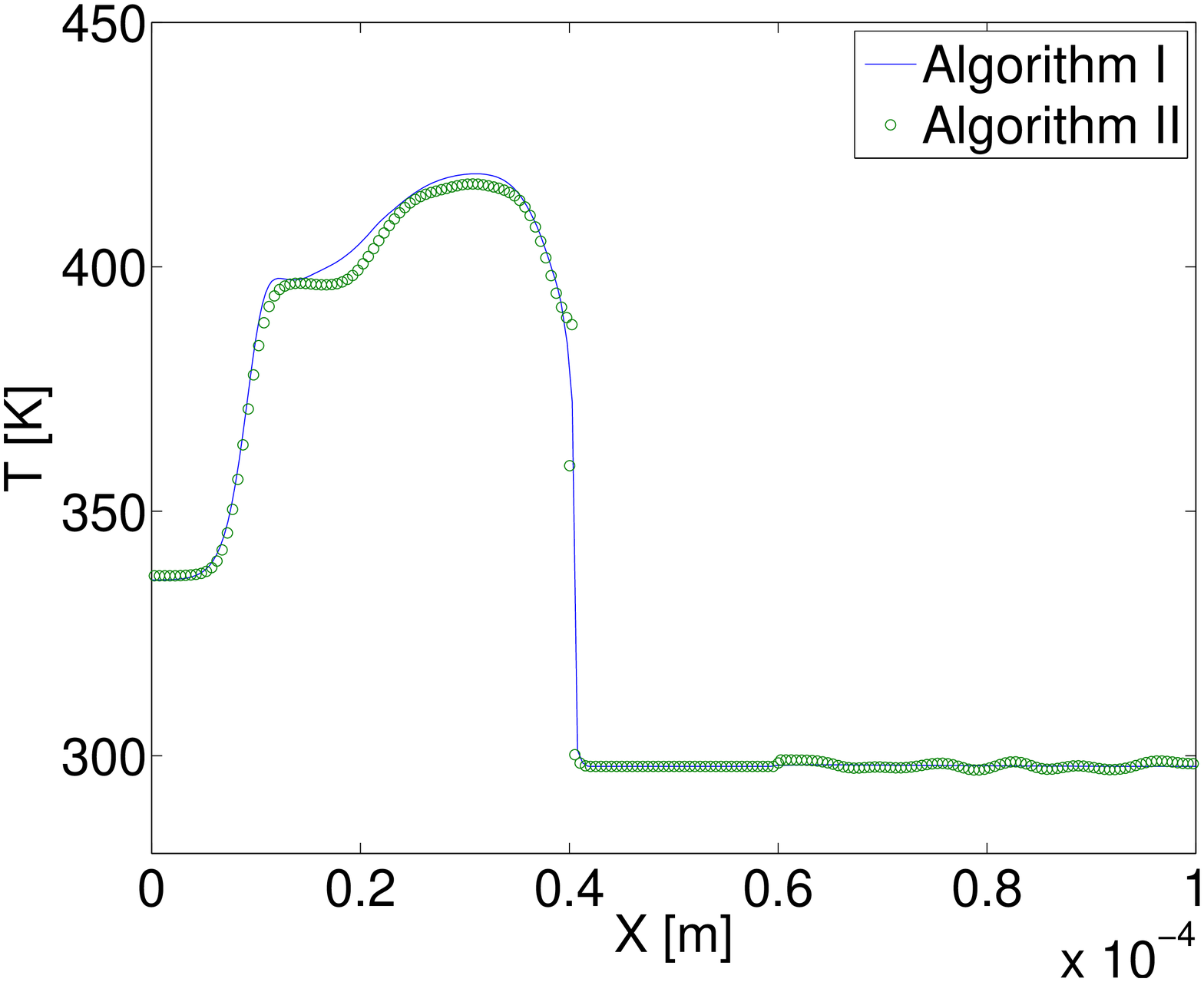,width=6.85cm}
 \end{center}
\caption{ \em Test 2: Logarithm of density(top left), pressure(top right), velocity(bottom left) and temperature(bottom right) for a liquid density $\rho_l = 1000~kg/m^3$ at time $t=1.75 \times 10^{-7}~s$. The solid lines are results from Algorithm I, the circles are from Algorithm II.}
\label{test2case1rho1to1000}
\end{figure}

\begin{figure}
\begin{center}
\epsfig{file= 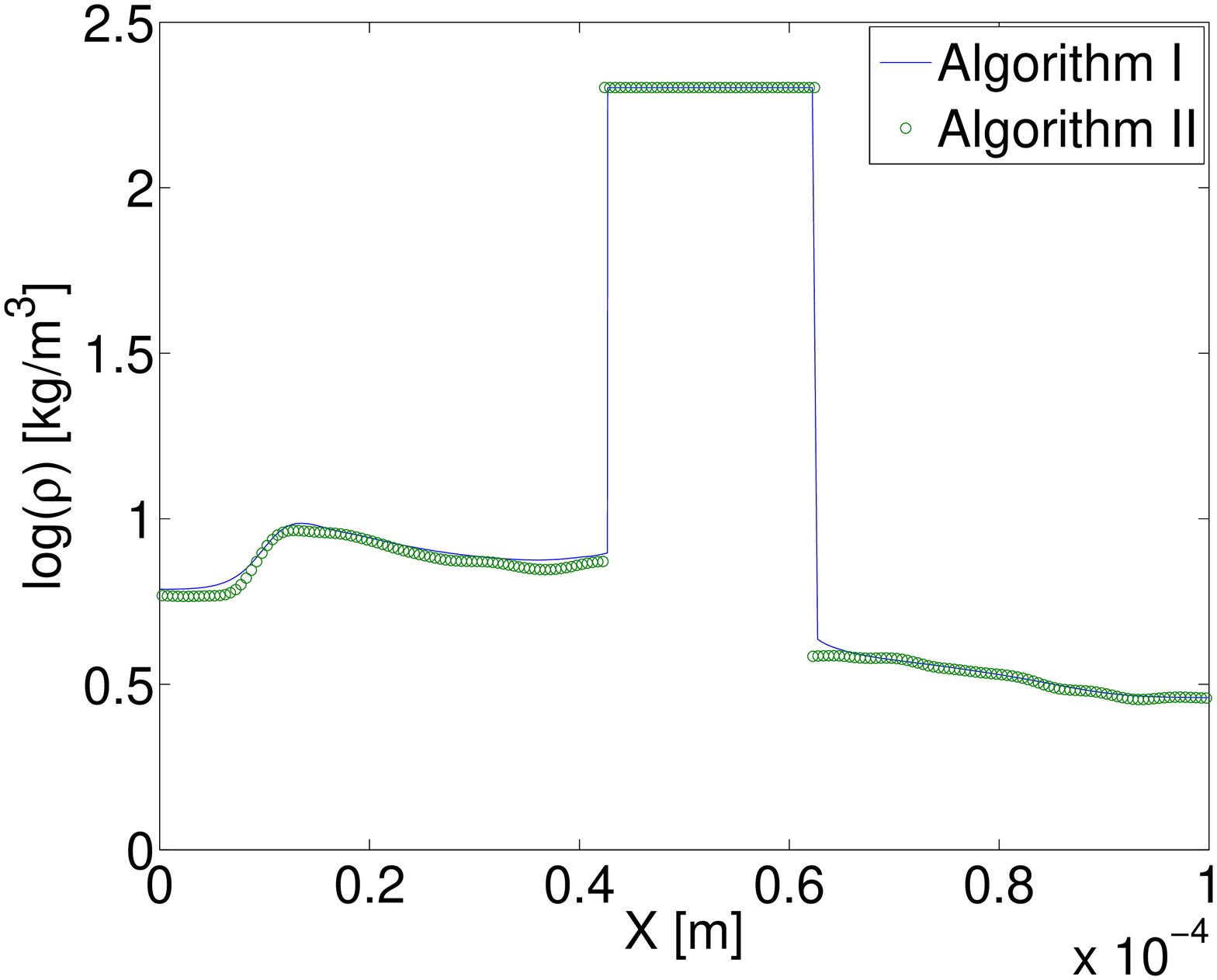,width=6.8cm}\hfill
\epsfig{file=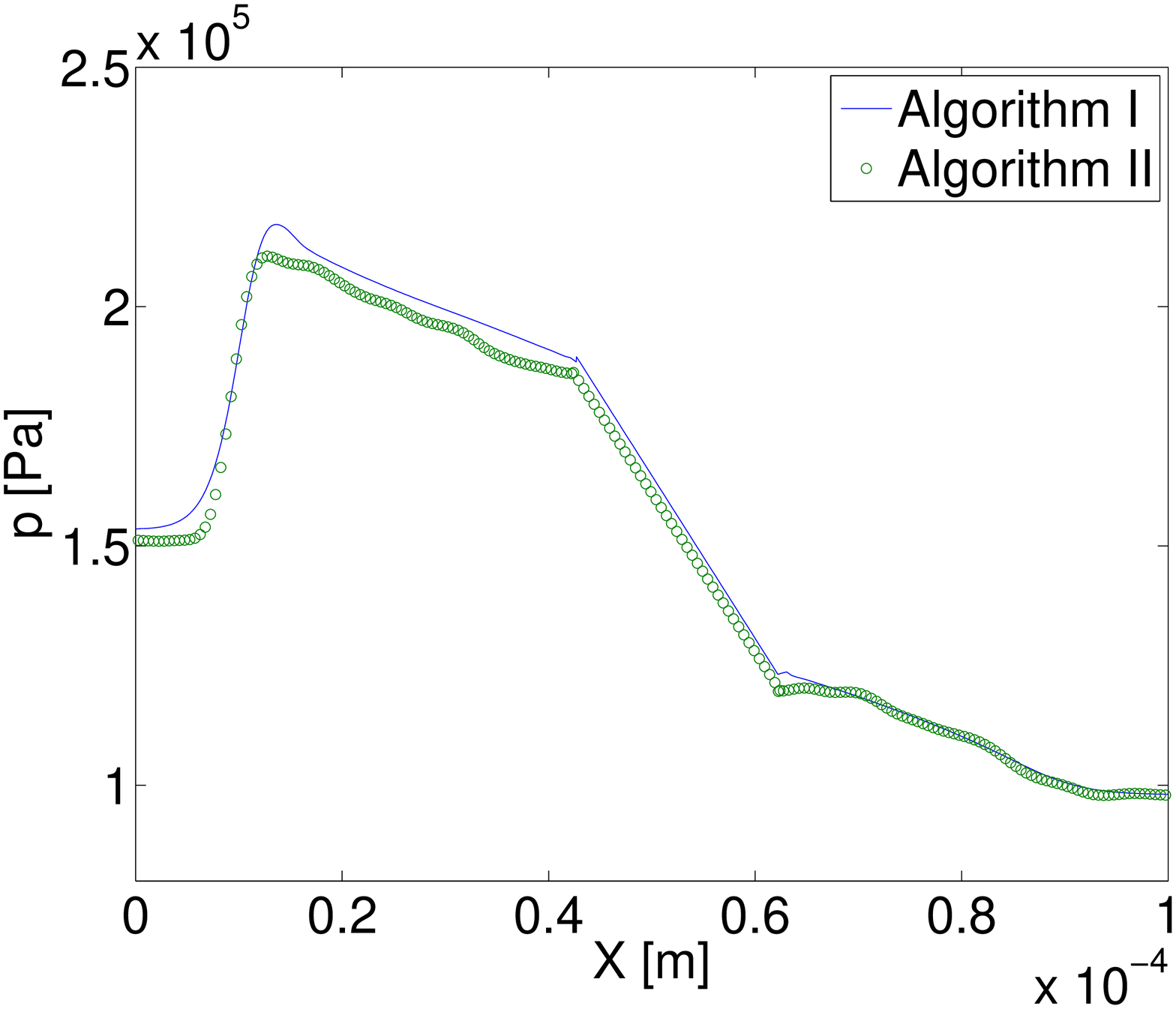,width=6.8cm}
\epsfig{file=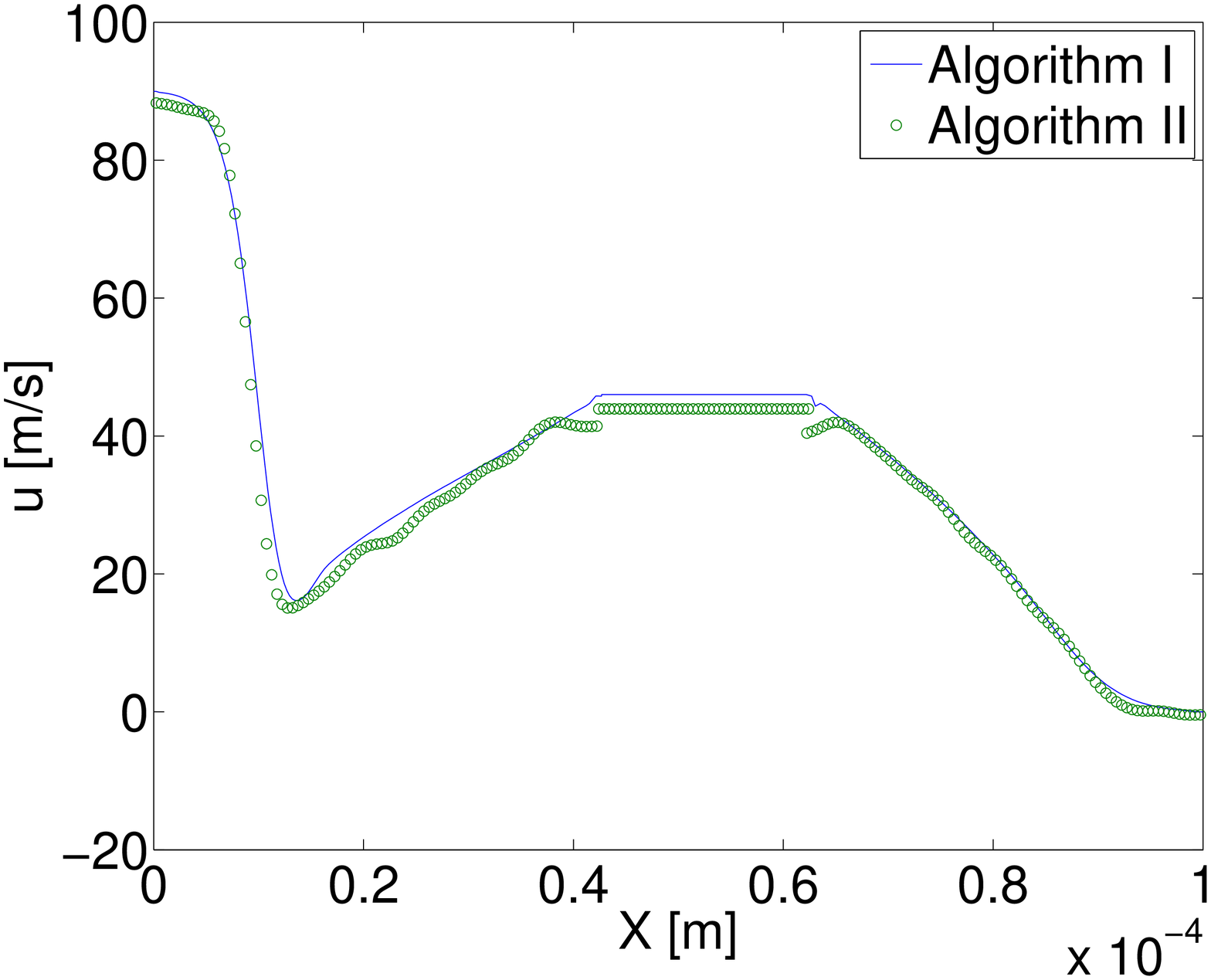,width=6.8cm}\hfill
\epsfig{file=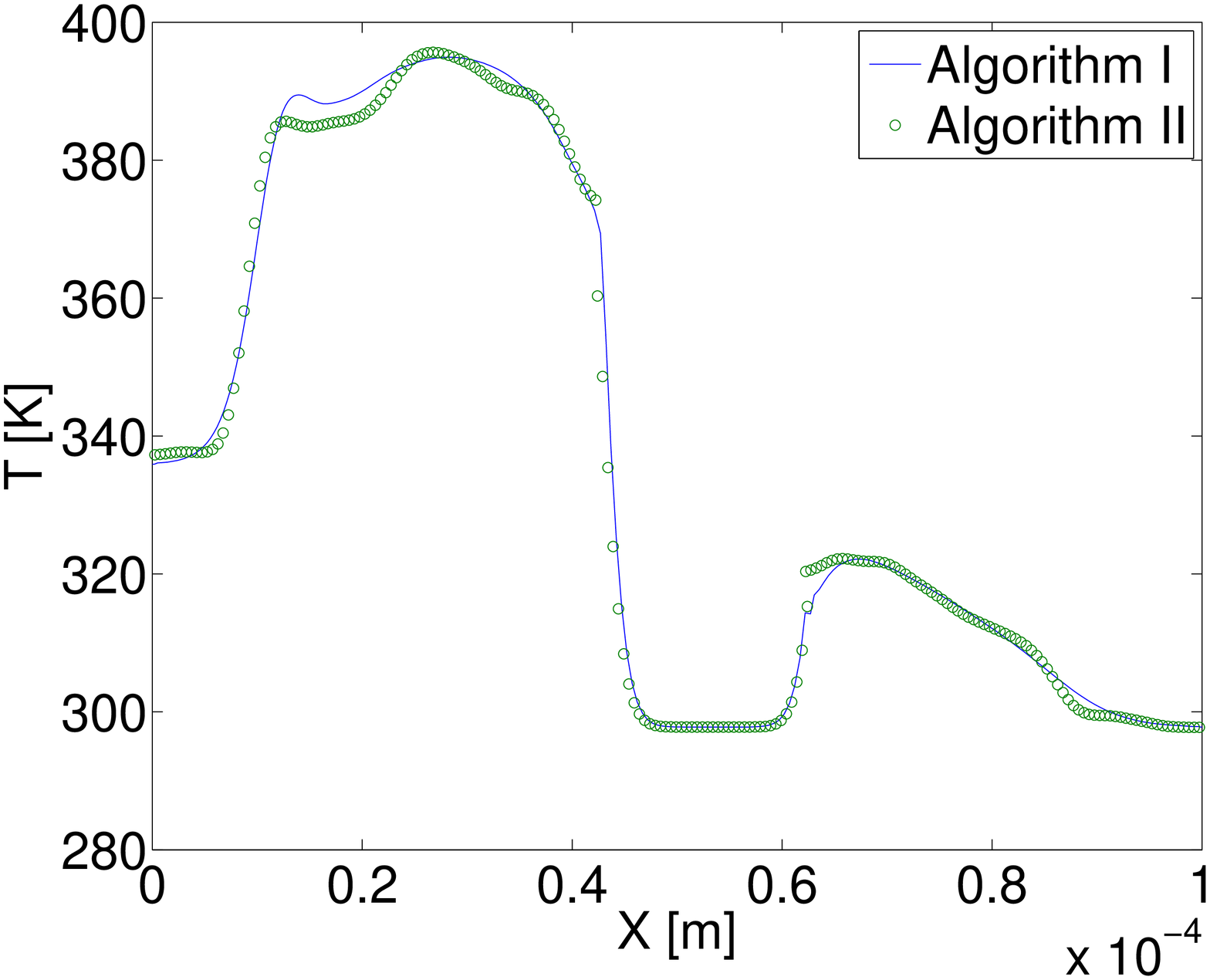,width=6.8cm}
 \end{center}
\caption{ \em Test 2: Logarithm of density(top left), pressure(top right), velocity(bottom left) and temperature(bottom right) for a liquid density $\rho_l = 10~kg/m^3$ at time $t=1.75 \times 10^{-7}~s$. The solid lines are results from Algorithm I, the circles are from Algorithm II.}
\label{test2case2rho1to10}
\end{figure}

For the Boltzmann equation we generate gas molecules in two ghost cells at the left boundary $x=a=0$ according 
to a Maxwellian with the initial parameters $\rho_l(0,a), {\bf u}_g(0,a) = 
(89.981,0,0)~m/s$ and $T_g(0,a)$. Similarly, we generate gas molecules in 
two cells at the right boundary $x=b$ according to a Maxwellian 
distribution with the parameters $\rho_g(0,b), T_g(0,b)$. The mean 
velocity ${\bf u}_g(t,b)$ is extrapolated from the neighboring cells at the left. 
After the free flow step, if gas molecules leave the boundaries 
at $a$ and $b$ we delete them. 
We note that we stopped the simulation before the wave reached the right boundary 
so that the mean velocity is still zero at this time. 
We have plotted the results at time 
$t=1.75\times 10^{-7}~s$ in figure \ref{test2case1rho1to1000}. 
One can observe that the shock wave travels to the right hitting the 
droplet, causing reflected as well as transmitted waves. 
The reflected wave is clearly visible in the pressure 
and the temperature plots of figure \ref{test2case1rho1to1000}, while both of the waves 
can be seen in figure \ref{test2case2rho1to10}. 
Again, the pressure profile is linear inside the droplet, however, with a reversed slope compared to the first test case. This originates from the fact that due to the interaction with the shock wave the accelaration of the droplet is now positive. The solutions obtained 
from both algorithms are very close to each other. 

In figure \ref{test2case2rho1to10} we have plotted the results from the same 
simulations with an liquid density $\rho_l = 10~kg/m^3$ at the same 
time $t=1.75 \times 10^{-7}~s$. It becomes apparent that for this lower liquid 
density the droplet has a lower inertia and is more easily displaced to the right. Compared to the initial state, 
the physical fields of the fluids have now been altered in almost the entire domain. 
Again, we see a 
good agreement of the solutions obtained from the coupling algorithms I 
and II.  

The above results are for the continuum gas regime. As expected, for  
the rarefied regime we see a discrepancy of the solutions obtained from 
algorithms I and II as in Test 1 for the largest Knudsen number. 
We do not explicitly present the corresponding results in this article. 

\subsection{Test 3}
In the final test case we again consider the same computational domain as in 
Test 2. The liquid droplet is initially at rest, 
occupying the domain $\Omega_l = [2\times 10^{-5}~m, 3\times 10^{-5}~m]$, while 
the rest of the domain is filled with gas. A similar test case has been 
studied in \cite{FMO98} for a unit interval.  
Initially, the gas is in thermal equilibrium with the initial state 
${\bf u}_g(0,x) = (0,0,0)~m/s, ~T_g(0,x) = 298~K$, and a spatially varying density $\rho_g(0,x)$. 
We assume $\rho_g(0,x) = 1~kg/m^3$ at the left of the droplet 
and a four times smaller value at the right of it. This means that the 
pressure at the left of the droplet is four times larger than at the right. 
The characteristic length is $1\times 10^{-5}~m$, and the Knudsen numbers 
at the left and right of the droplet are $0.011$ and $0.044$, respectively. 
The initial values for the liquid are $u_l(0,x) = 0m/s, ~T_l(0,x) = 298~K$, and 
the pressure is equal to the initial value in the gas at the right side of the droplet. 
We set the 
liquid density $\rho_l = 10~kg/m^3$.  
The boundary conditions are the same as in Test 1.
\begin{figure}
\begin{center}
\epsfig{file=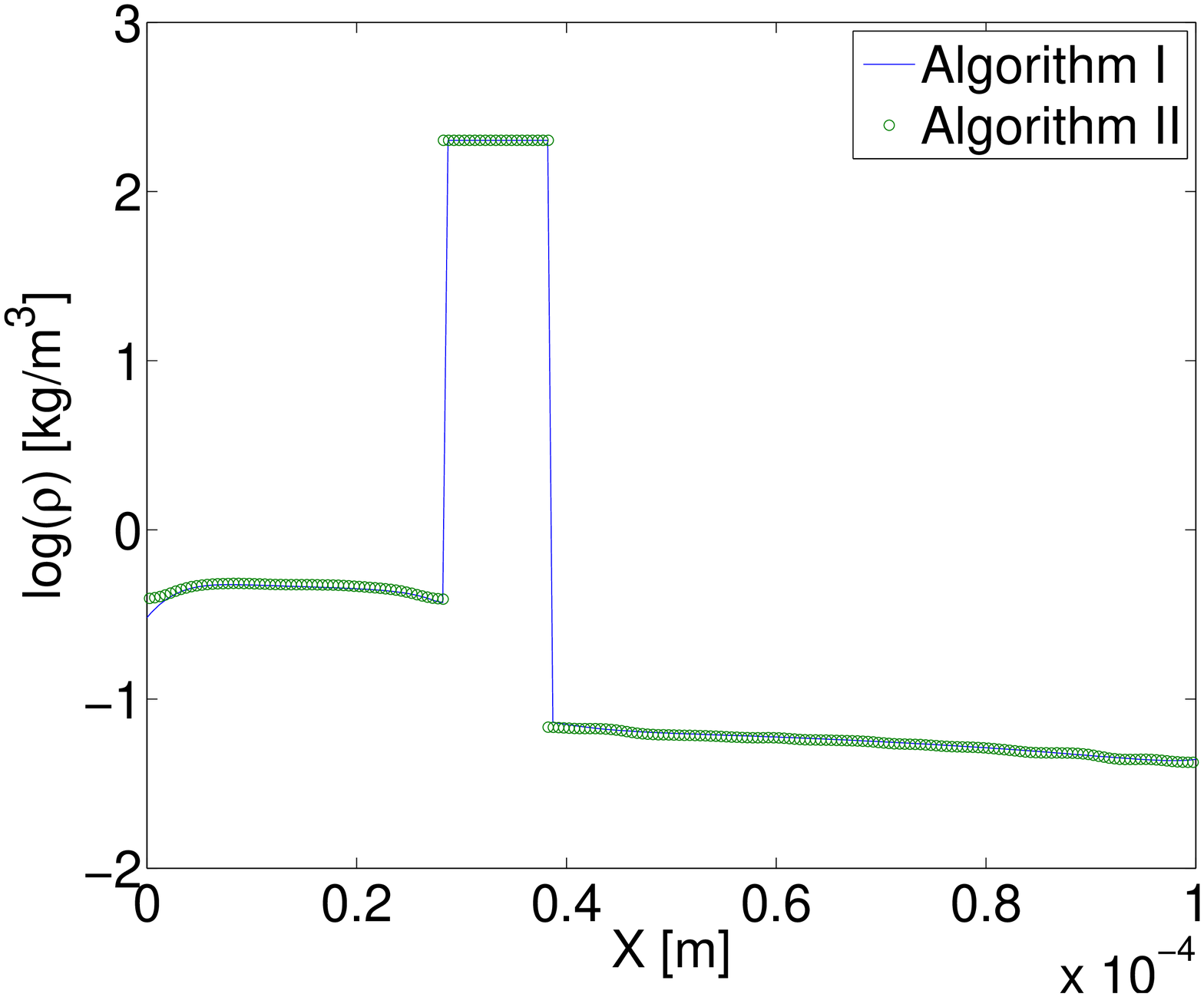,width=6.8cm}\hfill
\epsfig{file=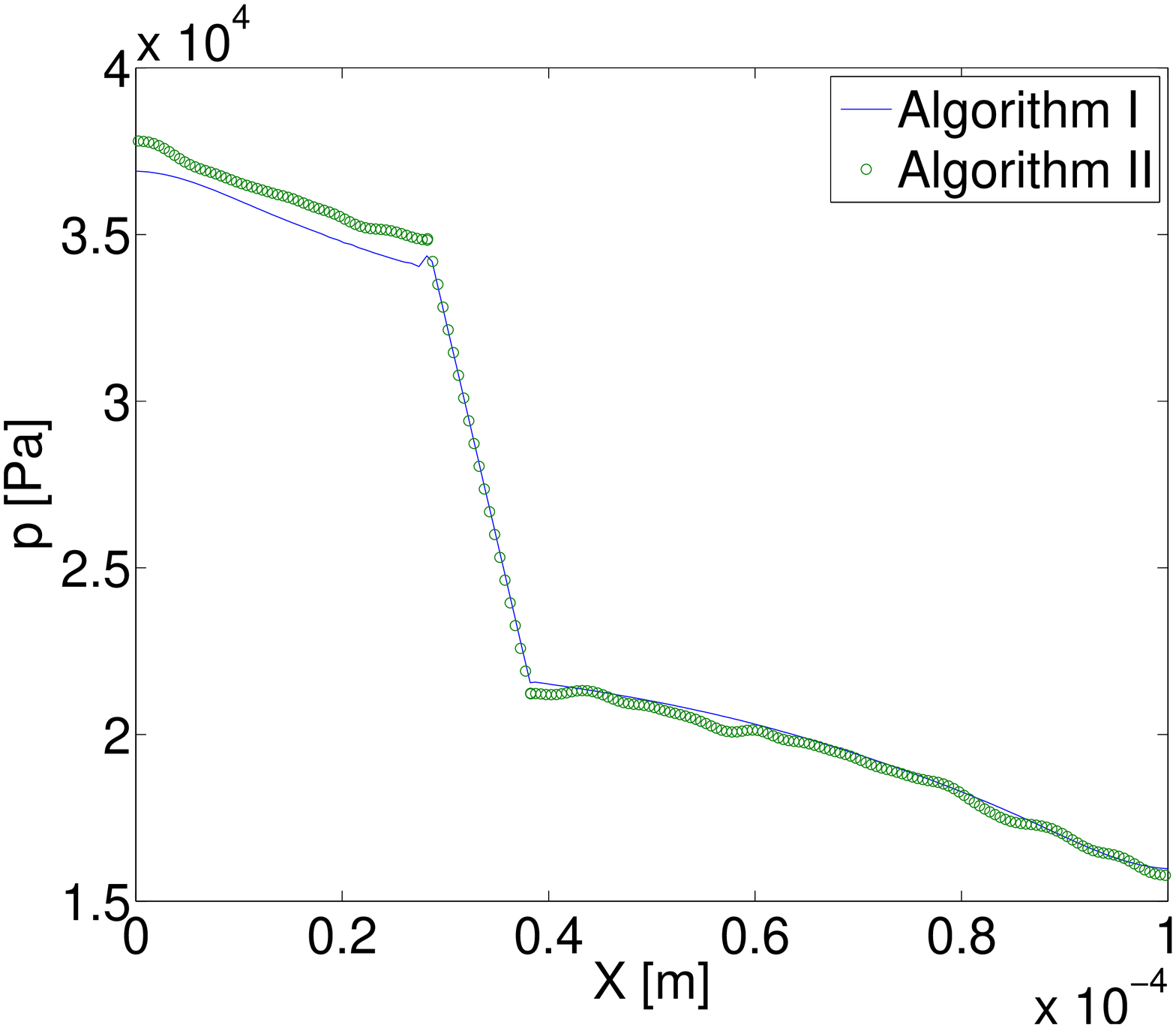,width=6.8cm}
\epsfig{file=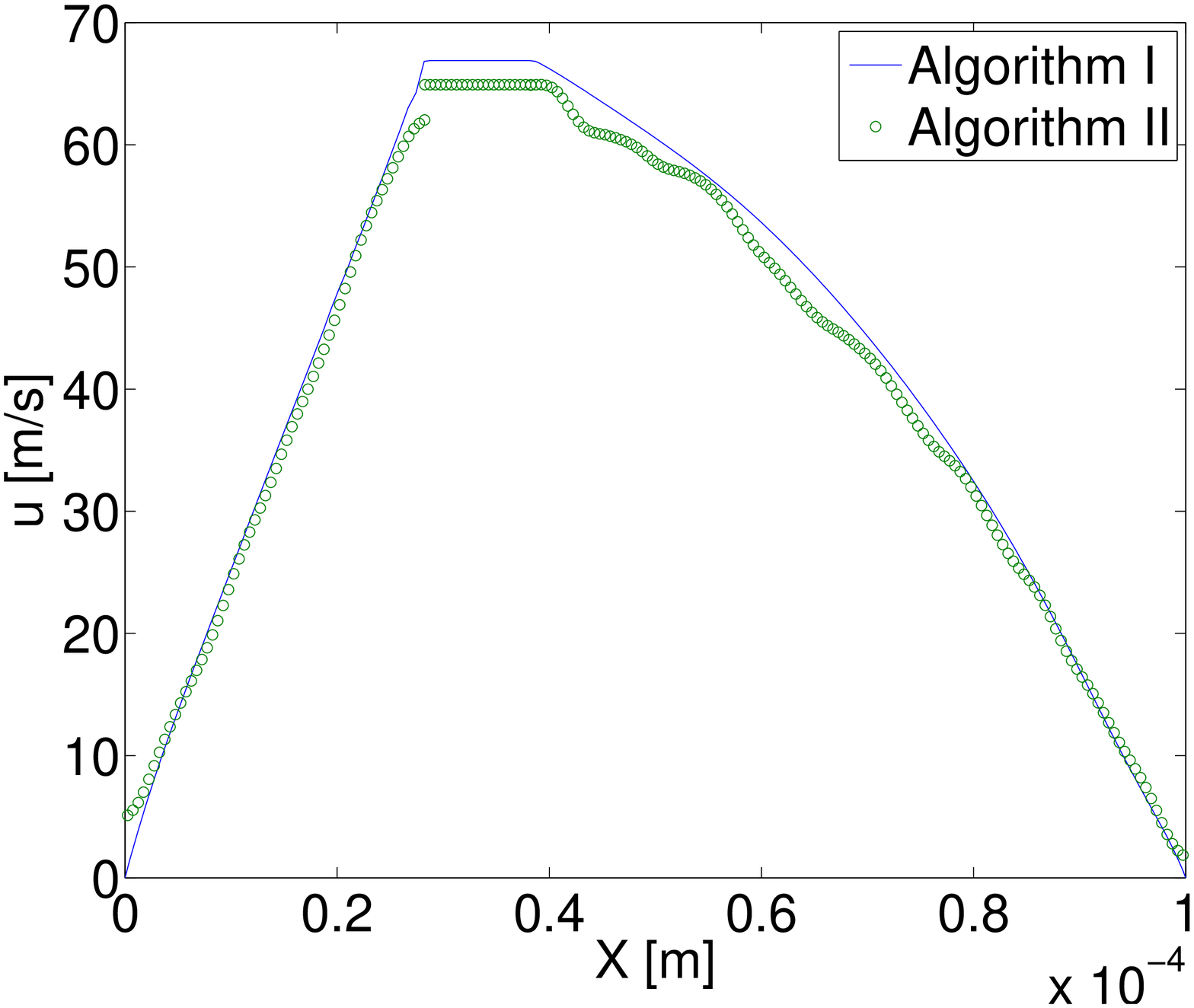,width=6.8cm}\hfill
\epsfig{file=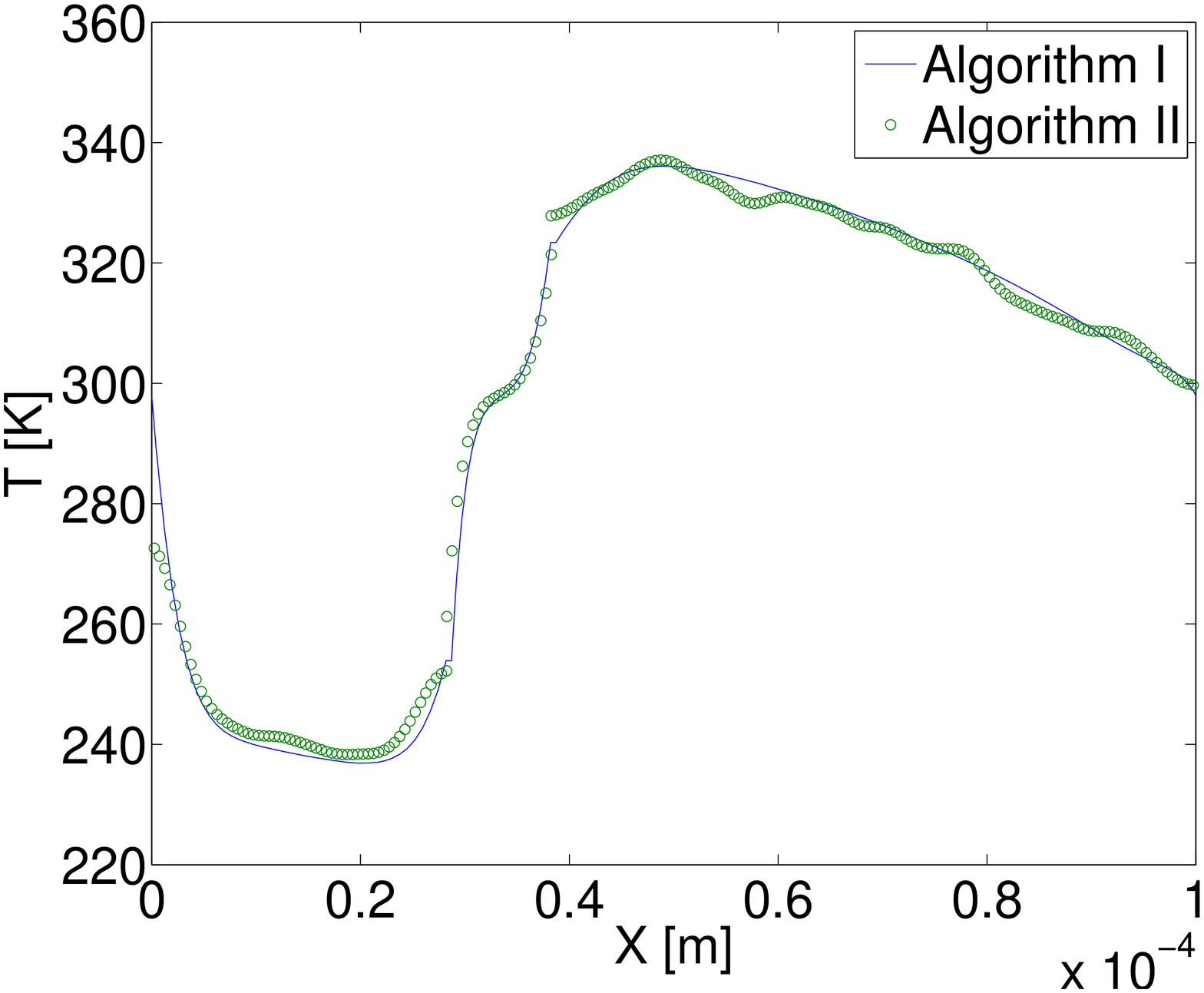,width=6.8cm}
 \end{center}
\caption{ \em Test 3: Logarithm of density(top left), pressure(top right), velocity(bottom left) and temperature(bottom right) for a liquid density $\rho_l = 10~kg/m^3$ at time 
$t=2.218\times 10^{-7}~s$. The solid lines are results from Algorithm I, the circles are from Algorithm II.  }
\label{test3case2rho1to10t2.218e-7}
\end{figure}
In figure \ref{test3case2rho1to10t2.218e-7} 
we have plotted the density, pressure, velocity and 
temperature after a short time $t=2.218\times 10^{-7}~s$. 
We see that the droplet has been slightly 
pushed to the right. The pressure on the left is higher. The 
gas has started to cool down on the left, and on the right it is getting compressed, 
resulting in a temperature increase.

\begin{figure}
\begin{center}
\epsfig{file=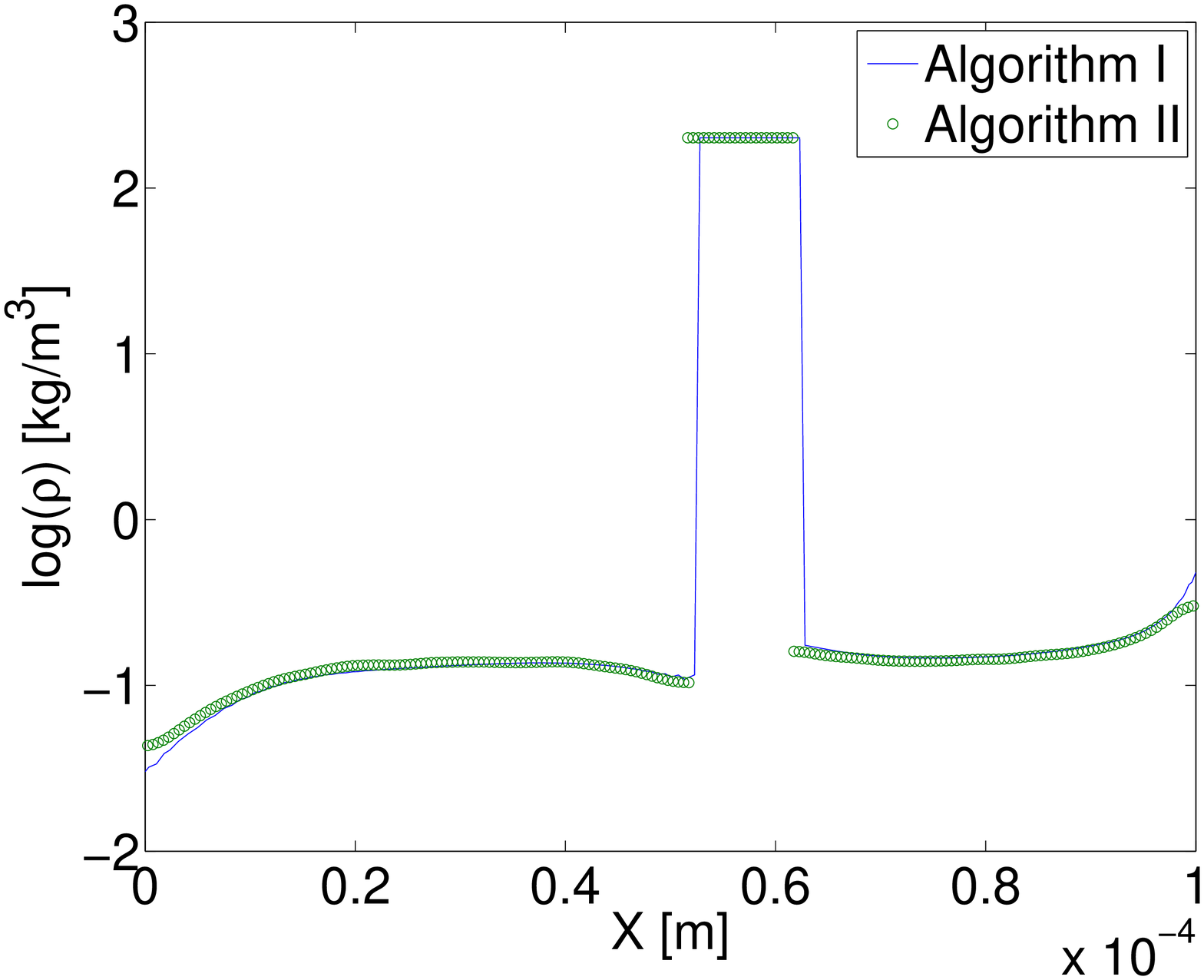,width=6.8cm}\hfill
\epsfig{file=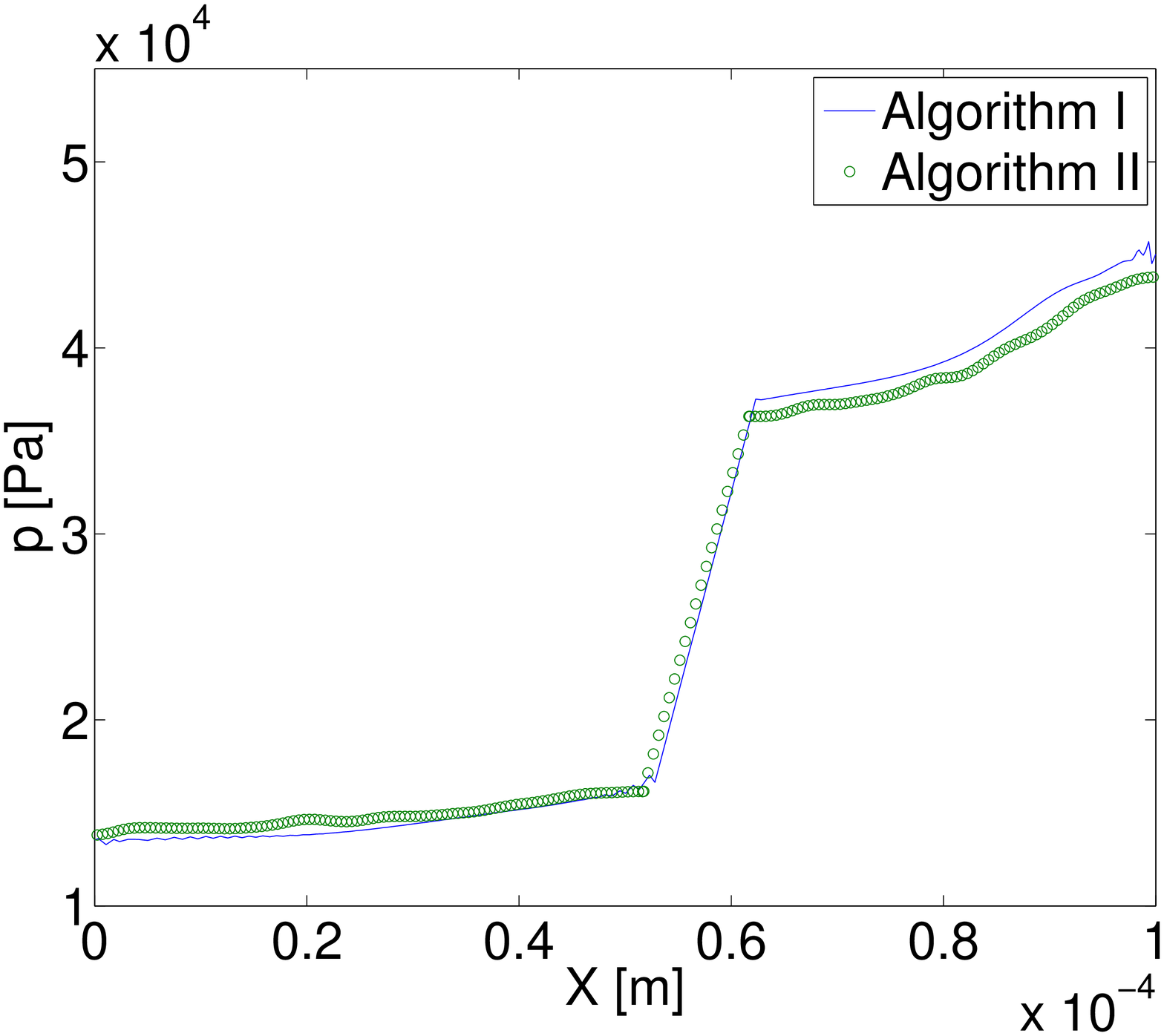,width=6.8cm}
\epsfig{file=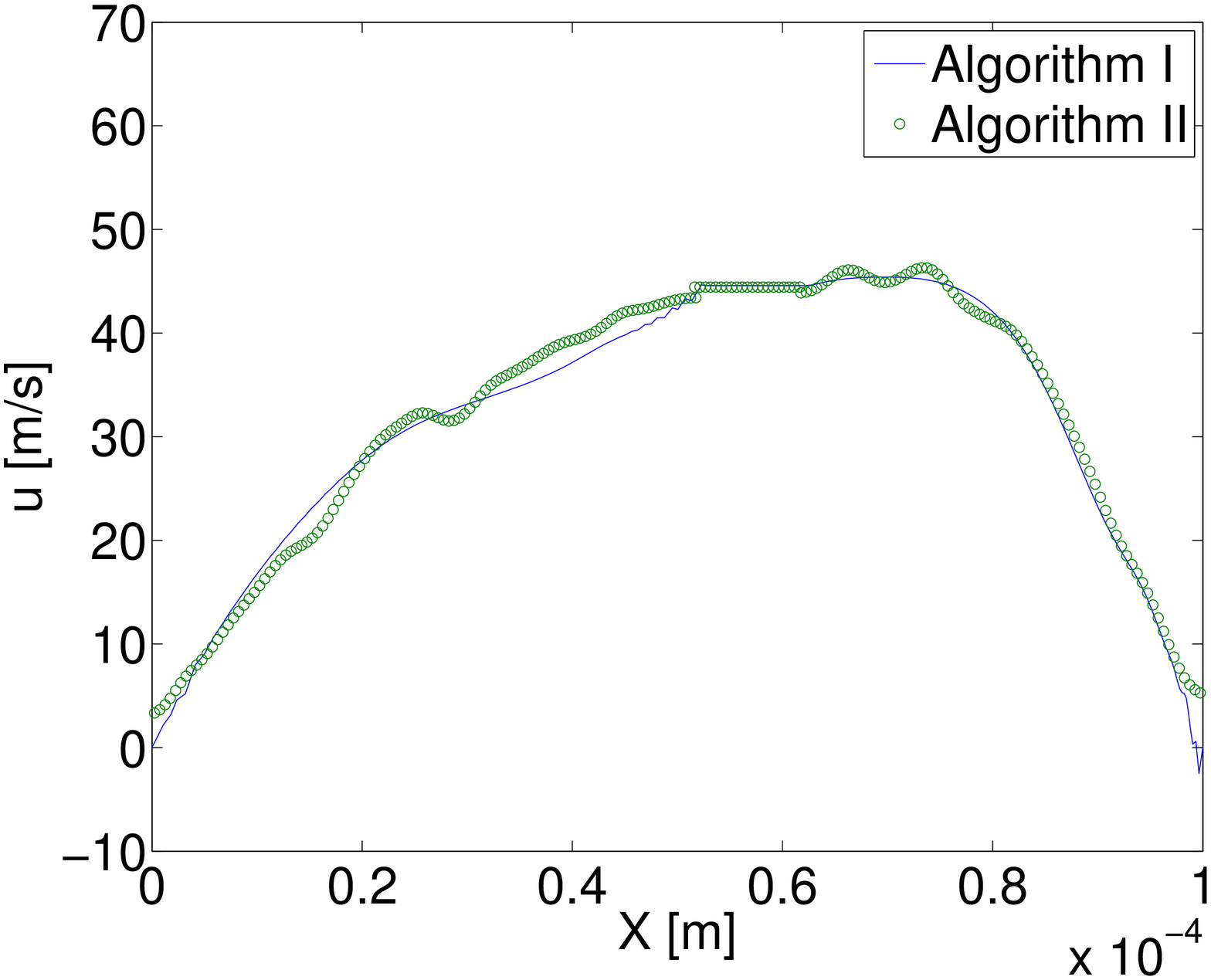,width=6.8cm}\hfill
\epsfig{file=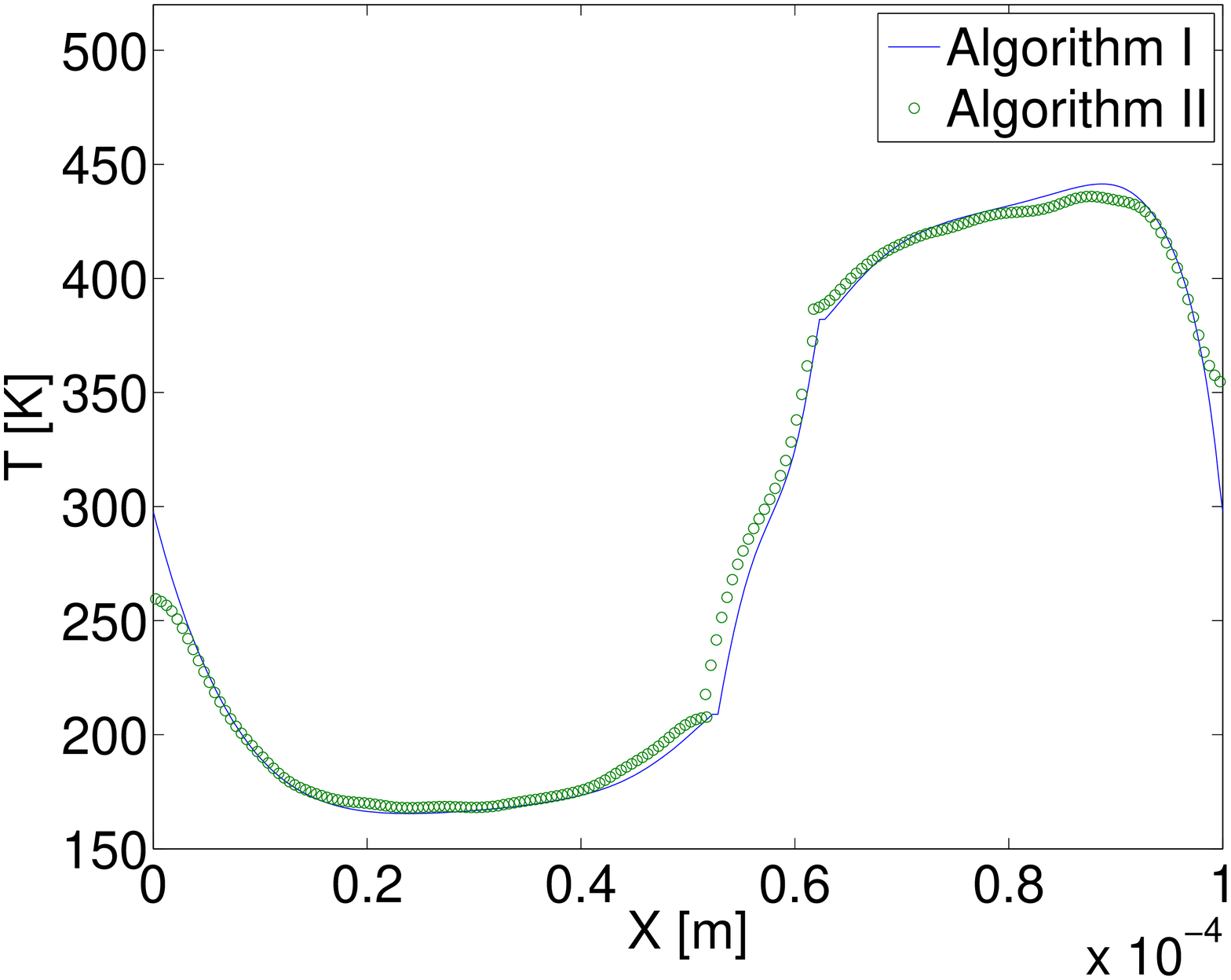,width=6.8cm}
 \end{center}
\caption{ \em Test 3: Logarithm of density(top left), pressure(top right), velocity(bottom left) and temperature(bottom right) for a liquid density $\rho_l = 10~kg/m^3$ at time 
$t=5.678\times 10^{-7}~s$. The solid lines are results from Algorithm I, the circles are from Algorithm II. }
\label{test3case2rho1to10t5.678e-7}
\end{figure}

At a later time $t=5.678\times 10^{-7}~s$ (figure 
\ref{test3case2rho1to10t5.678e-7}) 
we observe that the compression on the right has further increased and 
the gas has been further heated. 
Now the pressure on the right is higher 
causing a slowing down of the droplet. Going along with the transition from droplet accelaration to deceleration is a reversal of the slope
of the pressure profile inside the droplet. Later it will bounce back  
to the left. We observe that at both time levels the solutions obtained 
with both algorithms are in good agreement. 

\begin{figure}
\begin{center}
\epsfig{file=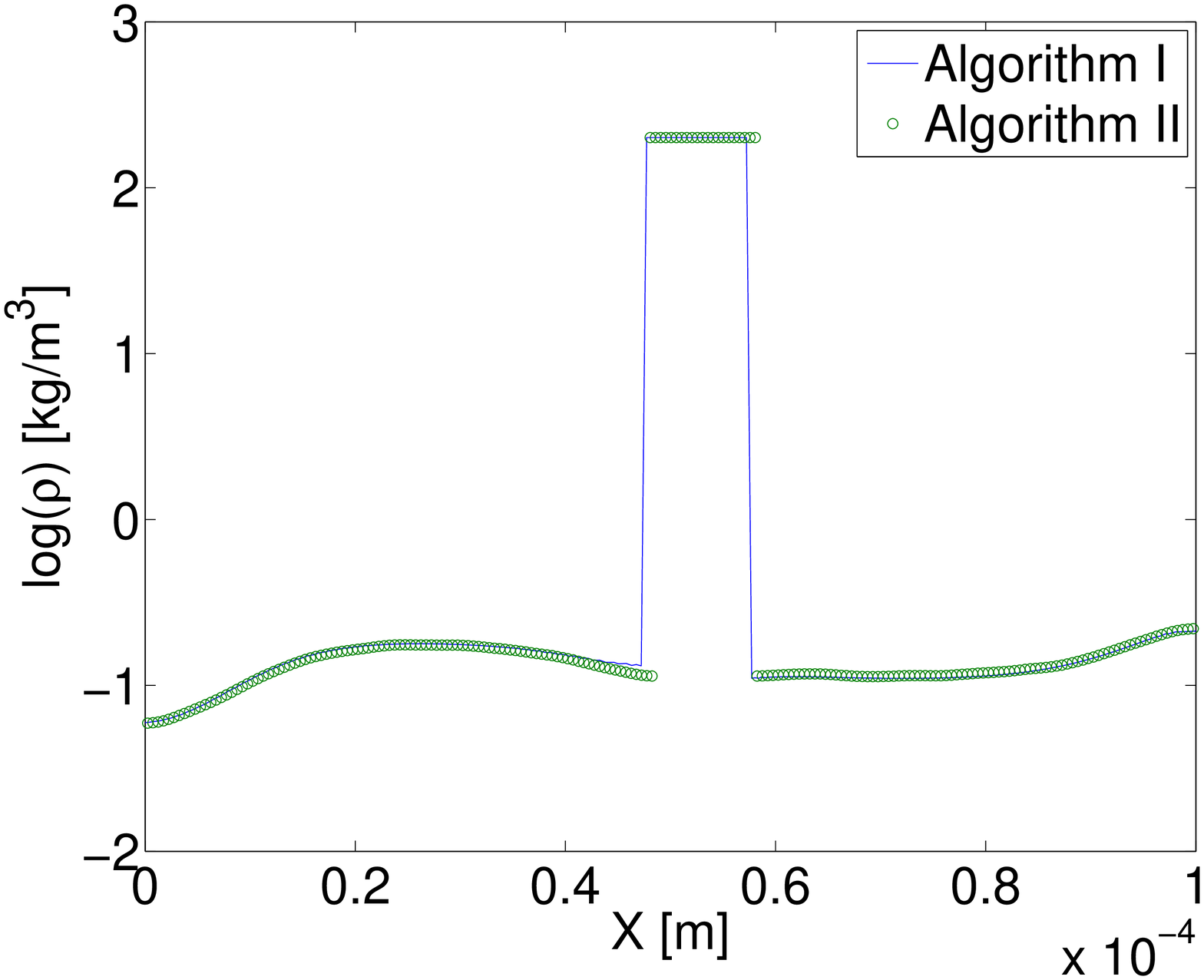,width=6.8cm}\hfill
\epsfig{file=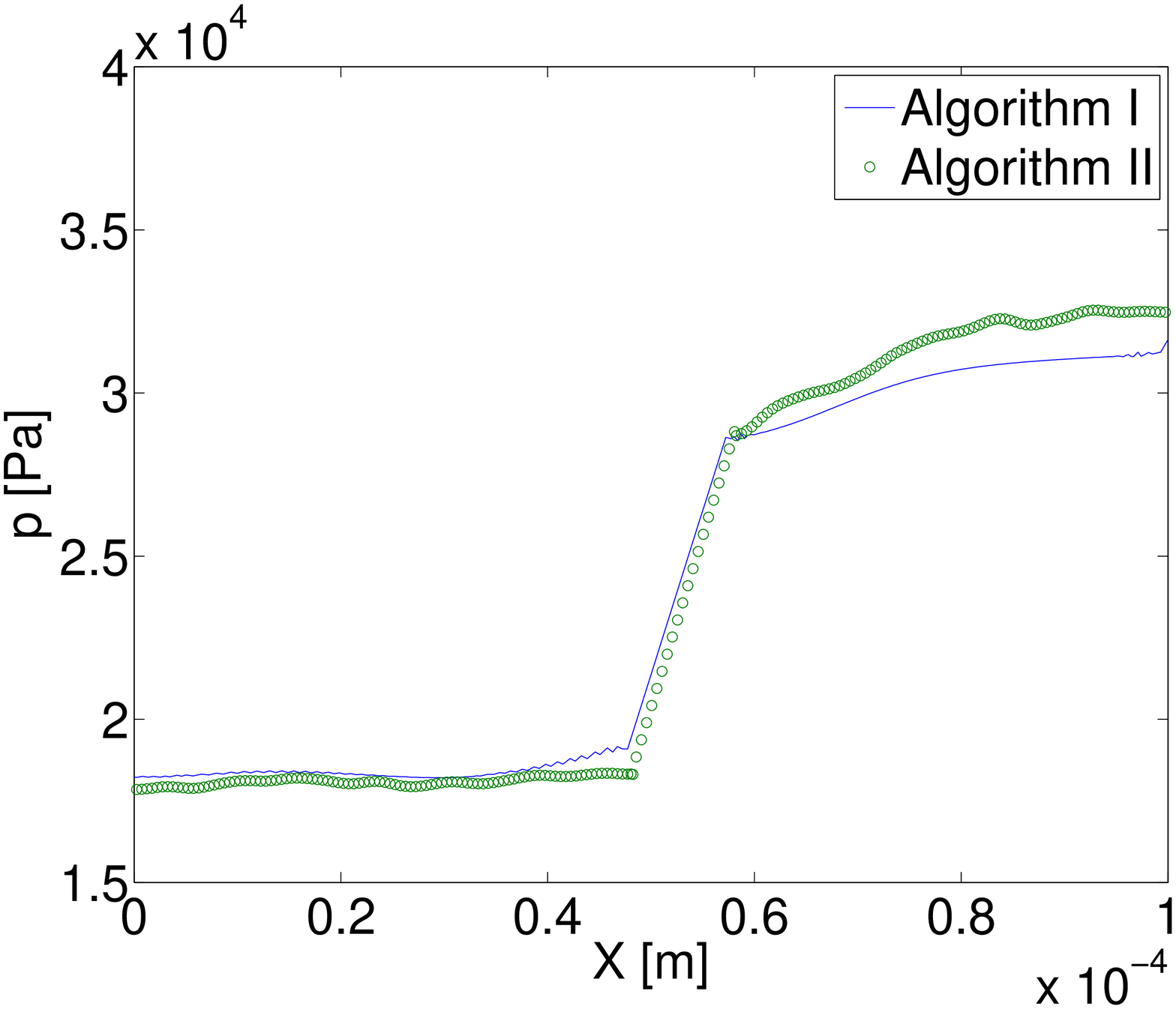,width=6.8cm}
\epsfig{file=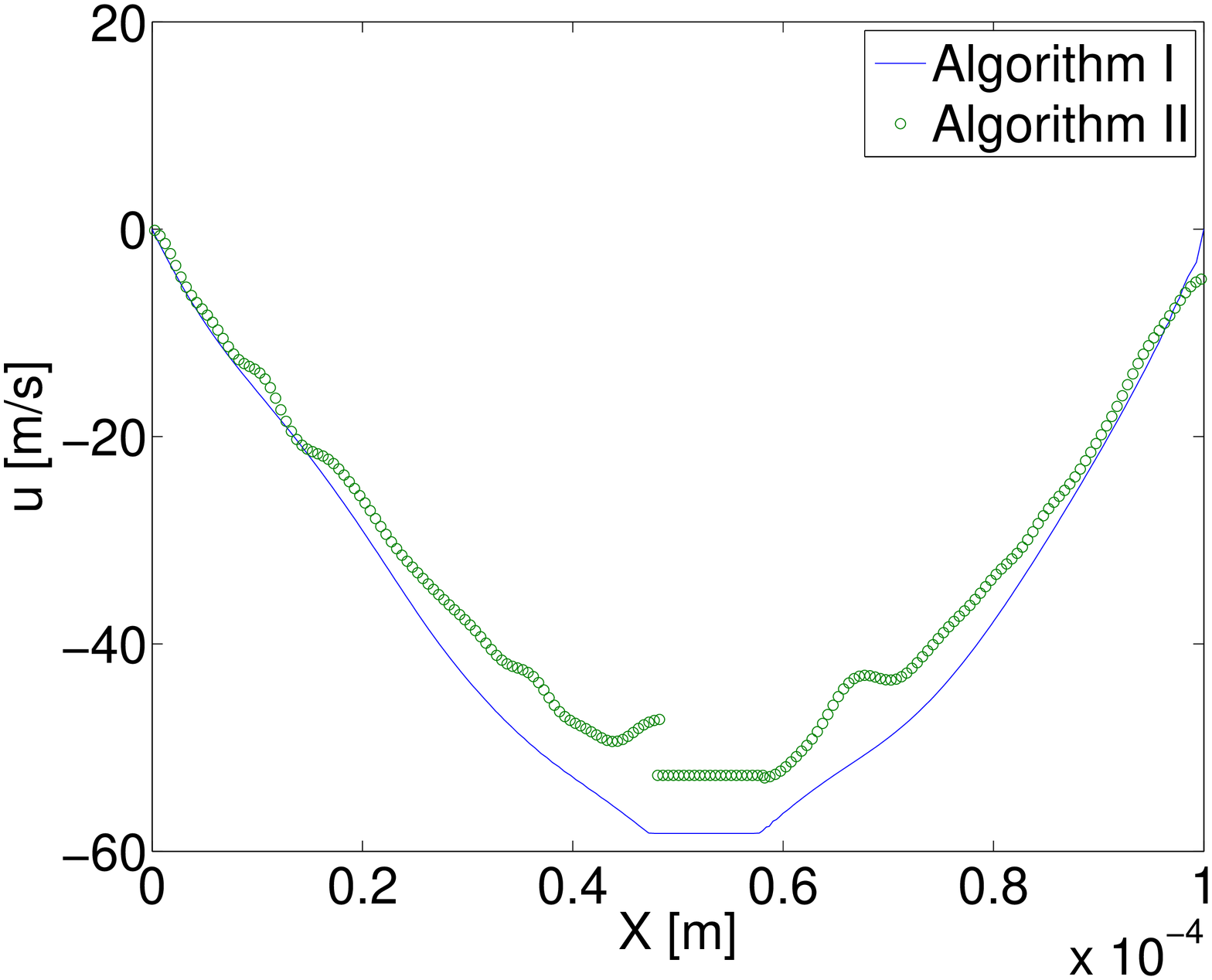,width=6.8cm}\hfill
\epsfig{file=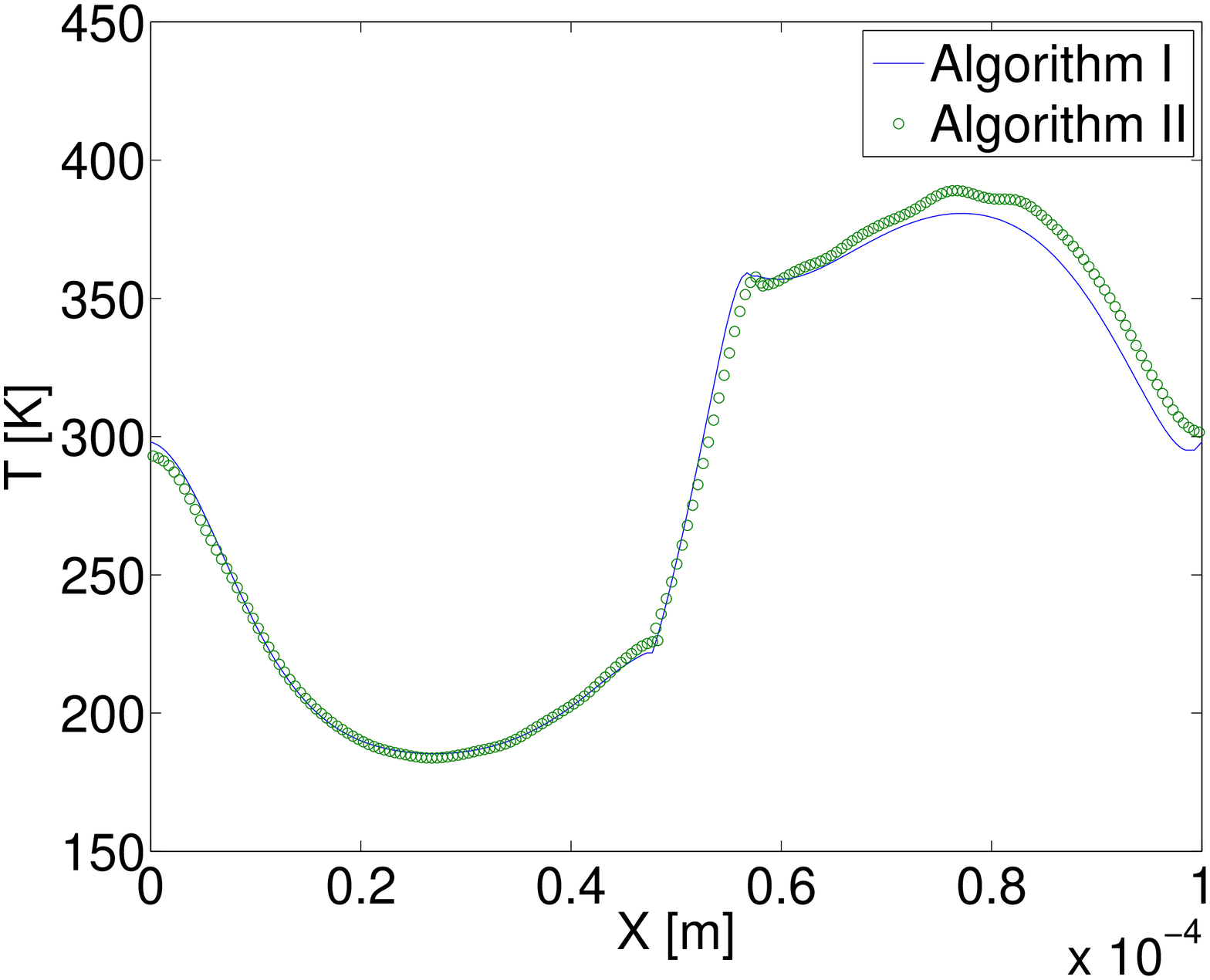,width=6.8cm}
 \end{center}
\caption{ \em Test 3: Logarithm of density(top left), pressure(top right), velocity(bottom left) and temperature(bottom right) for a liquid density $\rho_l = 10~kg/m^3$ at time 
$t=9.938\times 10^{-7}~s$. The solid lines are results from Algorithm I, the circles are from Algorithm II. }
\label{test3case2rho1to10t9.938e-7}
\end{figure}

At a still later time ($t=9.938\times 10^{-7}~s$) we see that the droplet is 
already pushed to the left, as displayed in 
figure \ref{test3case2rho1to10t9.938e-7}. 
At this time we see some discrepancy between algorithms I and II. This 
may be due to statistical fluctuation in the Boltzmann domain. 
However, the solutions are still comparable.  
If we compute further, the droplet will bounce back at the right boundary, and this 
oscillatory process will continue.

Similarly, for the rarefied gas regime, the solutions from Algorithms I and II 
do not match due to the breakdown of the continuum hypothesis inside the gas domain. 
As in Test 2 we do not present the results for larger 
Knudsen numbers.

\section{Conclusion and Outlook}

We have presented a coupling algorithm for a moving liquid phase inside a 
rarefied gas. The liquid phase is modeled by the incompressible Navier-Stokes 
equations, while the rarefied gas phase is modeled by the Boltzmann equation. 
The transport equations in the liquid are solved by a meshfree Lagrangian particle method, 
the Boltzmann equation by a DSMC type of particle method. Liquid 
particles overlap with the Boltzmann cells and the gas-liquid interface  
is determined from the particles at the boundary of the liquid domain. 
Interface boundary conditions are derived 
and a coupling algorithm for solving the Boltzmann equation
in combination with the incompressible 
Navier-Stokes equations is presented. To validate the coupling between the rarefied 
gas and the liquid we have modeled the gas by the compressible 
Navier-Stokes equations. 
Both the compressible and incompressible Navier-Stokes equations are 
solved by a meshfree Lagrangian particle method, called the 
Finite Pointset Method (FPM), where the gas and liquid particles 
are distinguished by assigning different flags.
We have also derived and implemented the interface 
boundary conditions for the compressible and incompressible Navier-Stokes 
equations which are well understood and widely used. 
In the continuum regime we have shown that the coupled solutions of the 
compressible and incompressible Navier-Stokes equations match well with 
the coupled solutions of the Boltzmann and the 
incompressible Navier-Stokes equations. 

Future work will concentrate on
the extension of the code to higher-dimensional problems, especially to 
the simulation of nano bubbles surrounded by an aqueos phase. Moreover, 
the method will be extended to simulate two-phase flows 
with evaporation and condensation \cite{BS04, HW08}.

\bibliographystyle{model1-num-names}
\bibliography{tiwarietal-BIBLIOGRAPHY}


\end{document}